\documentclass{amsart}

\usepackage{lmodern}
\usepackage{amsmath,amssymb,mathtools}

\usepackage{dsfont}
\usepackage[centertags]{amsmath}
\usepackage[usenames,dvipsnames]{xcolor}
\usepackage[colorlinks=true,linkcolor=Red,urlcolor=NavyBlue,citecolor=blue]{hyperref}
\usepackage{amsfonts}
\usepackage{MnSymbol}
\usepackage{amsthm}
\usepackage{newlfont}
\usepackage{amscd}
\usepackage{amsmath}
\usepackage{enumerate}
\usepackage[all,2cell]{xy}
\usepackage{tikz-cd}
\usetikzlibrary{positioning, calc}
\usepackage{mathrsfs}
\UseAllTwocells
\input xy
\xyoption{2cell}
\xyoption{all}
\usepackage{verbatim}
\usepackage{eucal}
\usepackage{graphicx}

\usepackage{pict2e}

\makeatletter
\newcommand{\adjunction}{\@ifstar\named@adjunction\normal@adjunction}
\newcommand{\normal@adjunction}[4]{%
	#1\colon #2%
	\mathrel{\vcenter{%
			\offinterlineskip\m@th
			\ialign{%
				\hfil$##$\hfil\cr
				\longrightharpoonup\cr
				\noalign{\kern-.3ex}
				\smallbot\cr
				\longleftharpoondown\cr
			}%
	}}%
	#3 \noloc #4%
}
\newcommand{\named@adjunction}[4]{%
	#2%
	\mathrel{\vcenter{%
			\offinterlineskip\m@th
			\ialign{%
				\hfil$##$\hfil\cr
				\scriptstyle#1\cr
				\noalign{\kern.1ex}
				\longrightharpoonup\cr
				\noalign{\kern-.3ex}
				\smallbot\cr
				\longleftharpoondown\cr
				\scriptstyle#4\cr
			}%
	}}%
	#3%
}
\newcommand{\longrightharpoonup}{\relbar\joinrel\rightharpoonup}
\newcommand{\longleftharpoondown}{\leftharpoondown\joinrel\relbar}
\newcommand\noloc{%
	\nobreak
	\mspace{6mu plus 1mu}
	{:}
	\nonscript\mkern-\thinmuskip
	\mathpunct{}
	\mspace{2mu}
}
\newcommand{\smallbot}{%
	\begingroup\setlength\unitlength{.15em}%
	\begin{picture}(1,1)
		\roundcap
		\polyline(0,0)(1,0)
		\polyline(0.5,0)(0.5,1)
	\end{picture}%
	\endgroup
}
\makeatother

\makeatletter
\newcommand*{\doublerightarrow}[2]{\mathrel{
		\settowidth{\@tempdima}{$\scriptstyle#1$}
		\settowidth{\@tempdimb}{$\scriptstyle#2$}
		\ifdim\@tempdimb>\@tempdima \@tempdima=\@tempdimb\fi
		\mathop{\vcenter{
				\offinterlineskip\ialign{\hbox to\dimexpr\@tempdima+1em{##}\cr
					\rightarrowfill\cr\noalign{\kern.5ex}
					\rightarrowfill\cr}}}\limits^{\!#1}_{\!#2}}}
\newcommand*{\triplerightarrow}[1]{\mathrel{
		\settowidth{\@tempdima}{$\scriptstyle#1$}
		\mathop{\vcenter{
				\offinterlineskip\ialign{\hbox to\dimexpr\@tempdima+1em{##}\cr
					\rightarrowfill\cr\noalign{\kern.5ex}
					\rightarrowfill\cr\noalign{\kern.5ex}
					\rightarrowfill\cr}}}\limits^{\!#1}}}
\makeatother

\newcommand{\NN}{\mathbb{N}}

\newcommand{\ZZ}{\mathbb{Z}}
\newcommand{\LL}{\mathbb{L}}

\def\sE{E}

\def\C{\mathcal{C}}

\def\F{\mathcal{F}}
\def\SS{\mathcal{S}}
\def\M{\mathcal{M}}
\def\N{\mathcal{N}}

\def\I{\mathcal{I}}
\def\cS{\mathcal{S}}
\def\T{\mathcal{T}}
\def\P{\mathcal{P}}
\def\Q{\mathcal{Q}}

\def\R{\mathcal{R}}
\renewcommand{\L}{\mathcal{L}}
\def\U{\mathcal{U}}

\def\g{\mathfrak{g}}

\def\lrarsimeq{\overset{\simeq}{\lrar}}

\def\l{\langle}
\def\r{\rangle}

\newcommand{\HSwarrow}{\kern0.05ex\vcenter{\hbox{\Huge\ensuremath{\Swarrow}}}\kern0.05ex}
\newcommand{\hSwarrow}{\kern0.05ex\vcenter{\hbox{\huge\ensuremath{\Swarrow}}}\kern0.05ex}
\newcommand{\LLSwarrow}{\kern0.05ex\vcenter{\hbox{\LARGE\ensuremath{\Swarrow}}}\kern0.05ex}
\newcommand{\LSwarrow}{\kern0.05ex\vcenter{\hbox{\Large\ensuremath{\Swarrow}}}\kern0.05ex}

\newcommand{\HSearrow}{\kern0.05ex\vcenter{\hbox{\Huge\ensuremath{\Searrow}}}\kern0.05ex}
\newcommand{\hSearrow}{\kern0.05ex\vcenter{\hbox{\huge\ensuremath{\Searrow}}}\kern0.05ex}
\newcommand{\LLSearrow}{\kern0.05ex\vcenter{\hbox{\LARGE\ensuremath{\Searrow}}}\kern0.05ex}
\newcommand{\LSearrow}{\kern0.05ex\vcenter{\hbox{\Large\ensuremath{\Searrow}}}\kern0.05ex}

\newcommand{\HDownarrow}{\kern0.05ex\vcenter{\hbox{\Huge\ensuremath{\Downarrow}}}\kern0.05ex}
\newcommand{\hDownarrow}{\kern0.05ex\vcenter{\hbox{\huge\ensuremath{\Downarrow}}}\kern0.05ex}
\newcommand{\LLDownarrow}{\kern0.05ex\vcenter{\hbox{\LARGE\ensuremath{\Downarrow}}}\kern0.05ex}
\newcommand{\LDownarrow}{\kern0.05ex\vcenter{\hbox{\Large\ensuremath{\Downarrow}}}\kern0.05ex}

\newcommand{\HUparrow}{\kern0.05ex\vcenter{\hbox{\Huge\ensuremath{\Uparrow}}}\kern0.05ex}
\newcommand{\hUparrow}{\kern0.05ex\vcenter{\hbox{\huge\ensuremath{\Uparrow}}}\kern0.05ex}
\newcommand{\LLUparrow}{\kern0.05ex\vcenter{\hbox{\LARGE\ensuremath{\Uparrow}}}\kern0.05ex}
\newcommand{\LUparrow}{\kern0.05ex\vcenter{\hbox{\Large\ensuremath{\Uparrow}}}\kern0.05ex}

\newcommand\restr[2]{{
		\left.\kern-\nulldelimiterspace 
		#1 
		\vphantom{\big|} 
		\right|_{#2} 
}}

\DeclareMathOperator{\hocolim}{hocolim}

\DeclareMathOperator{\Stab}{Stab}

\DeclareMathOperator{\Ext}{Ext}

\DeclareMathOperator{\Com}{Com}

\DeclareMathOperator{\Id}{Id}

\DeclareMathOperator{\id}{id}
\DeclareFontFamily{OT1}{pzc}{}
\DeclareFontShape{OT1}{pzc}{m}{it}{<-> s * [1.10] pzcmi7t}{}
\DeclareMathAlphabet{\mathpzc}{OT1}{pzc}{m}{it}
\DeclareMathOperator{\cof}{cof}

\DeclareMathOperator{\Alg}{Alg}
\DeclareMathOperator{\Ass}{Ass}
\DeclareMathOperator{\Free}{Free}
\DeclareMathOperator{\LMod}{LMod}

\DeclareMathOperator{\RMod}{RMod}

\DeclareMathOperator{\op}{op}

\DeclareMathOperator{\Map}{Map}

\DeclareMathOperator{\red}{red}

\DeclareMathOperator{\proj}{proj}

\DeclareMathOperator{\Sets}{Sets}

\DeclareMathOperator{\Cat}{Cat}

\DeclareMathOperator{\Mod}{Mod}

\DeclareMathOperator{\Fun}{Fun}

\DeclareMathOperator{\Ob}{Ob}
\DeclareMathOperator{\Sp}{Sp}
\DeclareMathOperator{\Ho}{Ho}
\DeclareMathOperator{\const}{const}

\DeclareMathOperator{\forget}{forgetful}

\DeclareMathOperator{\Op}{Op}

\DeclareMathOperator{\Ab}{Ab}

\DeclareMathOperator{\Tw}{Tw}
\DeclareMathOperator{\der}{h}
\DeclareMathOperator{\BMod}{BMod}

\DeclareMathOperator{\h}{h}

\DeclareMathOperator{\sMod}{sMod}

\DeclareMathOperator{\hocofib}{hocofib}

\DeclareMathOperator{\defi}{def}

\DeclareMathOperator{\IbMod}{IbMod}
\DeclareMathOperator{\ILMod}{I\ell Mod}
\DeclareMathOperator{\Hom}{Hom}
\DeclareMathOperator{\rL}{L}

\DeclareMathOperator{\sN}{N}

\DeclareMathOperator{\st}{t}

\DeclareMathOperator{\sB}{b}

\DeclareMathOperator{\Coll}{Coll}

\DeclareMathOperator{\Set}{Set}
\DeclareMathOperator{\sU}{U}

\DeclareMathOperator{\E}{E}

\DeclareMathOperator{\End}{End}
\DeclareMathOperator{\Fin}{Fin}

\DeclareMathOperator{\sS}{S}
\DeclareMathOperator{\HHH}{HH}
\DeclareMathOperator{\HHQ}{HQ}

\DeclareMathOperator{\ine}{in}
\DeclareMathOperator{\act}{act}
\DeclareMathOperator{\sur}{sur}
\DeclareMathOperator{\as}{as}
\DeclareMathOperator{\sr}{r}
\DeclareMathOperator{\Inf}{i}
\DeclareMathOperator{\me}{e}

\DeclareMathOperator{\nucom}{nCom}

\def\x{\overset}

\def\Stab{\textrm{Stab}}

\def\Lie{\textrm{Lie\,}}

\def\Hom{\textrm{Hom}}
\def\End{\textrm{End}}

\newcommand{\tgpd}{\kern0.05ex\vcenter{\hbox{\footnotesize\ensuremath{2}}}\kern0.05ex\mathcal{G}pd} 

\def\rar{\rightarrow}

\def\lrar{\longrightarrow}

\def\ovl{\overline}

\def\bar{\overline}

\newtheorem{theorem}{Theorem}[section]
\newtheorem{thm}{Theorem}[section]

\newtheorem{lem}[theorem]{Lemma}

\theoremstyle{definition}

\newtheorem{dfn}[theorem]{Definition}
\newtheorem{define}[theorem]{Definition}
\newtheorem{example}[theorem]{Example}

\newtheorem{examples}[theorem]{Examples}

\newtheorem{cons}[theorem]{Construction}
\newtheorem{conss}[theorem]{Constructions}

\theoremstyle{definition}

\newtheorem{rem}[theorem]{Remark}

\theoremstyle{definition}

\newtheorem{notn}[theorem]{Notation}

\theoremstyle{definition}

\newtheorem{notns}[theorem]{Notations}

\theoremstyle{corollary}

\newtheorem{cor}[theorem]{Corollary}

\theoremstyle{proposition}

\newtheorem{prop}[theorem]{Proposition}

\theoremstyle{observation}

\numberwithin{equation}{section}

\theoremstyle{definition}
\newtheorem{conv}[theorem]{Convention}

\theoremstyle{definition}

\begin{document}

\title[Hochschild and cotangent complexes of operadic algebras]%
{Hochschild and cotangent complexes \\ of operadic algebras}

\author{Truong Hoang}
\address{DEPARTMENT OF MATHEMATICS, HANOI FPT UNIVERSITY, VIET NAM.}

\email{truonghm@fe.edu.vn}



\begin{abstract} We make use of the cotangent complex formalism developed by Lurie to formulate Quillen cohomology of algebras over {an enriched operad}. Additionally, we introduce a spectral Hochschild cohomology theory for enriched operads and algebras over them. We prove that both the Quillen and Hochschild cohomologies of algebras over an operad can be controlled by the corresponding cohomologies of the operad itself. When passing to the category of simplicial sets, we assert that both these cohomology theories for operads, as well as their associated algebras, can be calculated in the same framework of spectrum valued functors on the twisted arrow $\infty$-category of the operad of interest.  Moreover, {we provide  {a convenient} cofiber sequence relating the Hochschild and cotangent complexes of an $\E_n$-space, establishing} an {unstable analogue} of a significant result obtained by Francis and Lurie. Our strategy introduces a novel perspective, focusing
	solely on the intrinsic properties of the operadic twisted arrow $\infty$-categories. 
	 
\end{abstract}

\maketitle

\tableofcontents

\section{Introduction}\label{s:introduction}

In \cite{Hoang}, we established the Quillen cohomology of operads, inspired by Harpaz-Nuiten-Prasma’s framework for enriched categories (\cite{YonatanCotangent}). In the present paper, we further extend these ideas to construct the Quillen cohomology of algebras over an operad, as well as a spectral version of Hochschild cohomology for such objects, with the aim of providing a broader framework for classical Hochschild cohomology.

\smallskip

\subsection{From tangent categories to spectral cotangent complexes}\label{s:tant}

Let us first recall the definitions of the classical \textbf{Hochschild cohomology} (\cite{Hoch}) and \textbf{Quillen cohomology} (\cite{Quillen}).  

\smallskip

Let $\textbf{k}$ be a commutative ring and let $A$ be a $\textbf{k}$-algebra. {Note} that the structure of a bimodule over $A$ is equivalent to that of a (left) module over $A\otimes A^{\op}$. Conceptually, the {$n$-th} \textit{Hochschild cohomology group} of $A$ with coefficients in an $A$-bimodule $N$ is given by
$$  \HHH^n(A;N) := \Ext^n_{A\otimes A^{\op}}(A,N).$$

\smallskip

Quillen cohomology is more universal. Let $X$ be an object in a model category $\M$. {We denote by} $\Ab(\M_{/X})$ the category of \textbf{abelian group objects} over $X${, and} assume that it comes equipped with a model structure transferred along the free-forgetful adjunction $\F :  \adjunction*{}{\M_{/X}}{\Ab(\M_{/X})}{} : \U$. Then one defines the \textbf{cotangent complex} of $X$ to be $\rL_X := \LL\F(X)$ the \textbf{derived abelianization} of $X$. Now the \textit{{$n$-th} Quillen cohomology group} of $X$ with coefficients in some $M \in \Ab(\M_{/X})$ is given by the formula
$$ \HHQ^n(X;M) := \pi_0 \Map^{\h}_{\Ab(\M_{/X})}(\rL_X , M[n]) $$
in which $M[n]$ signifies the $n$-suspension of $M$ in $\Ab(\M_{/X})$. 

\smallskip

Despite the successes, there are certain constraints within Quillen's theory (cf. \cite{YonatanCotangent, Hoang}). In response, Lurie contributes to the further advancement of Quillen's abelianization through the stabilization of $\infty$-categories. 

\smallskip

Let $X$ be an object living within a presentable $\infty$-category $\C$. The \textbf{tangent category to $\C$ at $X$} is by definition $\T_X\C := \Stab((\C_{/X})_*)$ the \textbf{stabilization} of $(\C_{/X})_*$, i.e. the pointed $\infty$-category associated to the over category $\C_{/X}$. Due to the presentability of $\C$, the canonical functor $\T_X\C \lrar \C_{/X}$ admits a left adjoint called the \textbf{suspension spectrum functor}, and denoted by $$\Sigma^{\infty}_+ : \C_{/X} \lrar \T_X\C.$$ In light of this, the \textbf{(spectral) cotangent complex} of $X$ is defined to be 
$$\rL_X := \Sigma^{\infty}_+(X).$$ Now the \textit{Quillen cohomology of} $X$ with coefficients in some $M \in \T_X\C$ is the space
$$ \HHQ^\star(X;M)  := \Map_{\T_X\C}(\rL_X , M),$$
and the {$n$\textit{-th}} \textit{Quillen cohomology group} is given by the formula $$\HHQ^{n}(X;M)  := \pi_0 \Map_{\T_X\C}(\rL_X , M[n]).$$ 
For more details about the above constructions, {we refer the reader} to \cite{Lurieha}. Additionally, according to \cite{YonatanBundle, YonatanCotangent}, one may opt to work within the framework of model categories.

\subsection{Cohomologies of enriched operads and algebras over them}\label{s:opalg}

Suppose we are given a sufficiently nice symmetric monoidal {model} category $\cS$. Let $\P$ be a $C$-colored operad enriched over $\cS${, where $C$ is a fixed set of colors}. As per Quillen-Lurie thinking, a cohomology theory for $\P$ depends on only the choice of a {model structure} on the category containing $\P$. What makes the situation intriguing is that $\mathcal{P}$ resides in various noteworthy categories.  {The following are of interest (see $\S$\ref{s:operad}, \ref{s:opmodules})}:

\smallskip

$\bullet$ $\P \in \Op(\cS)$: the category of $\cS$-enriched operads with non-fixed sets of colors,

\smallskip

$\bullet$ $\P \in \Op_C(\cS)$: the category of $C$-colored operads in $\cS$,

\smallskip

$\bullet$ ${\P^{\sB}} \in \BMod(\P)$:  the category of $\P$-bimodules, and

\smallskip

$\bullet$ ${\P^{\Inf}} \in \IbMod(\P)$:  the category of infinitesimal $\P$-bimodules.

\smallskip

{Here $\P^{\sB}$ (resp. $\P^{\Inf}$) denotes $\P$ itself considered as an object of $\BMod(\P)$ (resp. $\IbMod(\P)$).} In each case, an important invariant specific to $\mathcal{P}$ will emerge, as discussed below. 

\smallskip

We endow $\Op(\cS)$ with the \textbf{Dwyer-Kan {model structure}}, while the others are all equipped with the standard {model structure} transferred from that on $\cS$ {(cf. $\S$\ref{s:optransfermod})}. First, we refer to the Quillen cohomology of $\P\in\Op(\cS)$ as its \textbf{(proper) Quillen cohomology}, while the corresponding cohomology of $\P\in\Op_C(\cS)$ is called its \textbf{reduced Quillen cohomology}.

\smallskip

{Next, due to the fact that infinitesimal $\P$-bimodules are exactly the \textit{linearization} of the usual notion of $\P$-bimodules, it is natural to let the Hochschild complex of the operad $\P$ be classified by $\P^{\Inf} \in \IbMod(\P)$}. Nonetheless, note that the category $\IbMod(\P)$ remains unstable unless $\cS$ is stable. Therefore, to assure that the {resulting} cohomology theory is \textbf{abelian} (i.e. the cohomology groups are all abelian), we will let the cotangent complex of $\P^{\Inf}$, denoted by  $\rL_{\P^{\Inf}} \in \T_{\P^{\Inf}}\IbMod(\P)$, represent the \textbf{Hochschild complex of the operad} $\P$. More precisely, the Hochschild cohomology of $\P$ with coefficients in some $M\in \T_{\P^{\Inf}}\IbMod(\P)$ is the space
$$  \HHH^\star(\P;M)  := \Map^{\h}_{\T_{\P^{\Inf}}\IbMod(\P)}(\rL_{\P^{\Inf}} , M).$$
The above definition has full generality because when the category $\cS$ is stable then we obtain an equivalence $\T_{\P^{\Inf}}\IbMod(\P) \simeq \IbMod(\P)$, which identifies $\rL_{\P^{\Inf}}$ with ${\P^{\Inf}}$ itself.

\smallskip

Now let $A\in\Alg_\P(\cS)$ be a $\P$-algebra. Note that $A$ may also be regarded as an object of $\Mod_\P^A$ the category of $A$-\textbf{modules over} $\P$ {(see $\S$\ref{s:opmodules})}. {We will write $A^{\me}$ to denote $A$ itself in this role}. First, since we are mainly interested in the $\P$-algebra structure on $A$, it is obvious to exhibit the Quillen cohomology of $A\in\Alg_\P(\cS)$ as its \textbf{(proper) Quillen cohomology}. Moreover, according to the logic as in the above paragraph, we will exhibit $\rL_{A^{\me}} \in \T_{A^{\me}} \Mod_\P^A$, i.e. the cotangent complex of $A^{\me}$, as the \textbf{Hochschild complex of the} $\P$-\textbf{algebra} $A$. Now the Hochschild cohomology of $A$ with coefficients in some $N\in \T_{A^{\me}} \Mod_\P^A$ is the space
$$ \HHH^\star(A;N)  := \Map^{\h}_{\T_{A^{\me}} \Mod_\P^A}(\rL_{A^{\me}},N).$$

\subsection{Main statements}\label{s:intromainst}

We shall now summarize our main results, divided into three parts as follows.

\subsubsection{Formulations and comparison theorem}\label{s:formulaes} 

{It is noteworthy} that there exists a natural passage between infinitesimal $\P$-bimodules (resp. $\P$-bimodules) and $A$-modules over $\P$  (resp. $\P$-algebras) performed by the adjunction $(-) \circ_\P A \dashv \End_{A,-}$ (see $\S$\ref{s:endo}). Moreover, there is a commutative square of adjunctions between tangent categories:
\begin{equation}\label{eq:introkey}
	\begin{tikzcd}[row sep=3.5em, column sep =3.5em]
		\mathcal{T}_{\P^{\Inf}}\IbMod(\P) \arrow[r, shift left=1.5] \arrow[d, shift right=1.5, "(-)\circ^{\st}_\P A"'] & \mathcal{T}_{\P^{\sB}}\BMod(\P) \arrow[l, shift left=1.5, "\perp"'] \arrow[d, shift left=1.5, "(-)\circ^{\st}_\P A"]  & \\
		\T_{A^{\me}} \Mod_\P^A \arrow[r, shift right=1.5] \arrow[u, shift right=1.5, "\dashv", "\End_{A,-}^{\st}"'] & \T_A\Alg_\P(\cS)  \arrow[l, shift right=1.5, "\downvdash"] \arrow[u, shift left=1.5, "\vdash"', "\End_{A,-}^{\st}"] &
	\end{tikzcd}
\end{equation}
in which the vertical pairs are all induced by the adjunction $(-)\circ_\P A \dashv \End_{A,-}$ (see Notation \ref{no:circend}), and the horizontal pairs are given by the induction-restriction adjunctions. 

\smallskip

Under some mild assumptions on the triple $(\cS,\P,A)$, we obtain the following results.

\smallskip

{In Proposition \ref{p:bialg}, we establish the weak equivalences $$\rL_{\P^{\sB}}\circ^{\st}_\P A \simeq \rL_A \;\; \text{and} \;\; \rL_{\P^{\Inf}}\circ^{\st}_\P A \simeq \rL_{A^{\me}}$$
	in $\T_A\Alg_\P(\cS)$ and $\T_{A^{\me}}\Mod_\P^A$ respectively, where $\rL_{\P^{\sB}}$ denotes the cotangent complex of $\P$ as a bimodule over itself. These lead to the following \textbf{comparison theorem}.}

\begin{thm}(\ref{co:bialg}, \ref{t:opmain}) \label{co:introbialg} Suppose given a fibrant object $N\in \T_{A^{\me}} \Mod_\P^A$. 
	
	(i) There is a natural weak equivalence
	$$ \HHQ^\star(A ; N)  \simeq \Omega\HHQ^\star(\P ; \End^{\st}_{A,N})$$
	between the Quillen cohomology of $A$ with coefficients in $N$ and the loop space of Quillen cohomology of $\P$ with coefficients in $\End^{\st}_{A,N} \in \mathcal{T}_{\P^{\Inf}}\IbMod(\P)$.
	
	(ii) There is a natural weak equivalence
	$$ \HHH^\star(A ; N)  \simeq \HHH^\star(\P ; \End^{\st}_{A,N}) $$
	between the Hochschild cohomology of $A$ with coefficients in $N$ and Hochschild cohomology of $\P$ with coefficients in $\End^{\st}_{A,N}$. 
\end{thm}

\begin{rem} {When $\cS$ is stable, the cotangent complex $\rL_A \in \T_A\Alg_\P(\cS)$ can be represented as an $A$-module over $\P$. We then obtain a weak equivalence in $\Mod_\P^A$ of the form}
	{$ \ovl{\rL}_\P\circ_\P A \simeq \rL_A$}
	where $\ovl{\rL}_\P \in \IbMod(\P)$ is given at each $C$-sequence $(c_1,\cdots,c_m;c)$ by the formula
	$$\ovl{\rL}_\P(c_1,\cdots,c_m;c) = \P(c_1,\cdots,c_m;c)\otimes \hocolim_n\Omega^{n} [\, (\sS^{n})^{\otimes m}  \times_{1_\mathcal{S}}^{\h}  0\,].$$ 
	Here $\sS^{n}$ refers to the \textit{$n$-sphere} in $\cS$. {Moreover, we obtain weak equivalences, similar to those above:
	\begin{gather*} \HHQ^\star(A ; N)  \simeq \Omega\HHQ^\star(\P ; \End_{A,N}),  \;\; \text{and} \\ \HHH^\star(A ; N)  \simeq \HHH^\star(\P ; \End_{A,N})	
	\end{gather*}	
		with the only difference being that $N$ is now an object of $\Mod_\P^A$.} 
\end{rem}

\subsubsection{Simplicial operads and algebras over them}\label{s:simp} 

We will now assume that $\cS$ is the category $\Set_\Delta$ of simplicial sets, and that $\P$ is a fibrant and $\Sigma$-cofibrant $C$-colored operad in $\Set_\Delta${, for some fixed set of colors $C$}. Recall from \cite{Hoang} that there is a chain of equivalences of $\infty$-categories
\begin{gather*}
	\T_\P\Op(\Set_\Delta)_\infty \simeq \T_\P\Op_C(\Set_\Delta)_\infty \simeq \T_{\P^{\sB}}\BMod(\P)_\infty \\
	\simeq \T_{\P^{\Inf}}\IbMod(\P)_\infty \simeq \Fun(\Tw(\P) , \Sp)
\end{gather*}
where $\Tw(\P)$ signifies the \textbf{twisted arrow $\infty$}-\textbf{category} of $\P$ and $\Sp$ denotes the $\infty$-\textbf{category of spectra}. 

\smallskip

Due to the above, both Quillen and Hochschild cohomologies of $\P$ can be calculated in the framework of functors $\Tw(\P) \lrar \Sp$. In {\cite{Hoang}}, we have proved that $\rL_\P \in  \T_\P\Op(\Set_\Delta)_\infty$ is identified with $\F_\P[-1] \in \Fun(\Tw(\P) , \Sp)$ where the functor $\F_\P : \Tw(\P) \lrar \Sp$ takes any operation $\mu\in\P$ of arity $m$ to $\F_\P(\mu) = \mathbb{S}^{\oplus m}${, and while} the Hochschild complex $\rL_{\P^{\Inf}} \in \T_{\P^{\Inf}}\IbMod(\P)_\infty$ is identified with $\ovl{\mathbb{S}}: \Tw(\P) \lrar \Sp$ the constant functor on the \textit{sphere spectrum} (see Proposition \ref{p:hoch}). {We further obtain the following}.

\begin{thm}\label{t:introsimpal}(\ref{t:endan}) {Suppose that $\P \in \Op_C(\Set_\Delta)$ is fibrant and $\Sigma$-cofibrant. Let $A$ be a cofibrant $\P$-algebra and let $N \in \T_{A^{\me}} \Mod_\P^A$ be a fibrant object.}
	
	\smallskip

	(i) The Quillen cohomology of $A$ with coefficients in $N$ is computed by
	$$ \HHQ^\star(A ; N) \simeq \Omega\HHQ^\star(\P ; \End^{\st}_{A,N}) \simeq \Map_{\Fun(\Tw(\P) , \Sp)}(\F_\P, \widetilde{\End}^{\st}_{A,N})$$
	{where $\widetilde{\End}^{\st}_{A,N}: \Tw(\P) \lrar \Sp$ is the functor naturally associated with the pair $(A,N)$ (see Construction \ref{cons:endan})}.
	
	\smallskip
	
	(ii) The Hochschild cohomology of $A$ with coefficients in $N$ is computed by
	$$ \HHH^\star(A ; N) \simeq \HHH^\star(\P ; \End^{\st}_{A,N}) \simeq \Map_{\Fun(\Tw(\P) , \Sp)}(\ovl{\mathbb{S}}, \widetilde{\End}^{\st}_{A,N}).$$
	In particular, the $n$-th Hochschild cohomology group is given by 
	$$\HHH^n(A ; N) \cong \pi_{-n} \, \lim\widetilde{\End}^{\st}_{A,N}.$$
\end{thm}

\subsubsection{Quillen principle for $\E_n$-spaces}\label{ss:qprin} While Quillen cohomology plays a central role in the study of deformation theory and obstruction theory, its actual computation is far from straightforward. In addressing this complexity, one approach is to establish a connection between Quillen cohomology and the notably more accessible Hochschild cohomology. We will say
that a given object of interest is subject to the \textbf{Quillen principle} if there exists a convenient (co)fiber
sequence relating the cotangent and Hochschild complexes of that object. 

\smallskip

Let $A$ be an $\E_n$-algebra in $\Set_\Delta$. We let $\Free_A(*)$ denote the free $A$-module over $\E_n$ generated by a singleton, and let $\eta_A : \Free_A(*) \lrar A^{\me}$ be the map in $\Mod_{\E_n}^A$ classified by the unit of $A$. We have a canonical map $\Sigma^\infty_+(\eta_A) \lrar \rL_{A^{\me}}$ where $\Sigma^\infty_+(\eta_A)$ refers to the image of $\eta_A$ through the functor $$\Sigma^\infty_+ : (\Mod_{\E_n}^A)_{/A^{\me}} \lrar \T_{A^{\me}}\Mod_{\E_n}^A.$$

\begin{thm}\label{t:introQprin3} (\ref{t:Qprin3}, \ref{co:Qprin3})  {Suppose that $A$ is a cofibrant $\E_n$-algebra. Then} there is a cofiber sequence in $\T_{A^{\me}}\Mod_{\E_n}^A$ of the form
	\begin{equation}\label{eq:introQprin3}
		\Sigma^\infty_+(\eta_A) \lrar \rL_{A^{\me}} \lrar \rL_A[n]
	\end{equation}
	where we use the same notation for the derived image of $\rL_A\in\T_{A}\Alg_{\E_n}(\Set_\Delta)$ under the equivalence $\T_{A}\Alg_{\E_n}(\Set_\Delta) \simeq \T_{A^{\me}}\Mod_{\E_n}^A$. Consequently, for a given fibrant object $\N \in \T_{A^{\me}}\Mod_{\E_n}^A$, we obtain a fiber sequence connecting the two cohomologies 
	$$ \Omega^n\HHQ^\star(A;\N) \lrar \HHH^\star(A;\N) \lrar |\N|$$
	in which $|\N| := \{1_A\} \times^{\h}_A \N_{0,0}$, i.e. the homotopy fiber over $1_A$ of the structure map $\N_{0,0} \lrar A^{\me}$.
\end{thm}

This assertion  {provides an unstable analogue} of a result obtained by Francis and Lurie, which establishes the Quillen principle for $\E_n$-algebras in stable categories, such as chain complexes or symmetric spectra (cf. {\cite[Theorem 2.26]{Francis}} and {\cite[Theorem 7.3.5.1]{Lurieha}}).

\bigskip

{\textbf{\underline{Acknowledgements}.} The author would like to thank the referees for their careful reading of the manuscript and for their insightful comments, which have greatly contributed to improving the clarity and presentation of the paper.}

\section{Background and conventions}\label{s:background}

This section is devoted to some preliminaries relevant to the theory of enriched  operads and their associated modules, as well as the standard homotopy theories of such objects. We also revisit the tangent category formalism, from which we may obtain the central concepts of the paper including the spectral cotangent complex and Quillen cohomology.

\subsection{Enriched operads}\label{s:operad}

{Let $(\SS, \otimes, 1_\cS)$ be a \textbf{symmetric monoidal category} whose underlying category $\cS$ is cocomplete and whose monoidal product $-\otimes-$ distributes over colimits. We will refer to $\SS$ as the \textbf{base category}.}

\smallskip

{Let $C\in \Sets$ be a given set, regarded as the \textbf{set of colors}}. By definition, a \textbf{symmetric} $C$-\textbf{collection} in $\SS$ consists of a collection $$M = \{M(c_1,\cdots,c_n;c) \, | \, c_i,c \in C , n\geqslant0  \}$$ of objects in $\cS$ such that for each sequence $(c_1,\cdots,c_n;c)$ and each permutation $\sigma\in \Sigma_n$, there is a map of the form $$\sigma^{*} : M(c_1,\cdots,c_n;c) \lrar M(c_{\sigma(1)},\cdots,c_{\sigma(n)};c) .$$ Such maps {are required to establish a \textit{right $\Sigma_n$-action}, meaning that for every $\sigma,\tau \in \Sigma_n$, one has $\sigma^{*}\tau^{*} = (\tau\sigma)^*$ and $\iota_n^* = \Id$ where $\iota_n \in \Sigma_n$ denotes  the trivial permutation}. We let $\Coll_C(\SS)$ denote  the\textbf{ category of symmetric $C$-collections in $\SS$}.

\begin{rem} When $C$ is a singleton then a symmetric $C$-collection will be called a $\Sigma_*$-\textbf{object}. More concretely, a $\Sigma_*$-object $M$ consists of a collection $\{M(n)\}_{n\geq0}$ such that each $M(n)$ is a right $\Sigma_n$-object in $\cS$. The \textbf{category of $\Sigma_*$-objects in} $\cS$ will be denoted by $\Sigma_*(\cS)$.
\end{rem}

The category $\Coll_C(\SS)$ admits a monoidal structure given by the \textbf{composite product}
$$ -\circ- : \Coll_C(\SS) \times \Coll_C(\SS) \lrar \Coll_C(\SS)$$
{(see, e.g., \cite{Yonatan})}. We will denote by $\I_C$ the monoidal unit which agrees with $\emptyset_\cS$ (i.e. the initial object of $\cS$) on all the levels except that $\I_C(c;c)=1_{\cS}$ for every $c\in C$.

\begin{dfn} A \textbf{symmetric} $C$-\textbf{colored operad} in $\cS$ is a monoid in the monoidal category $(\Coll_C(\SS),-\circ-,\I_C)$. We will write $\Op_C(\SS)$ to denote the \textbf{category of symmetric $C$-colored operads in} $\cS$. When $C$ is a singleton, we will write $\Op_*(\SS)$ {for} the corresponding category. 
\end{dfn}

\begin{conv} {All operads considered in this paper are assumed to be symmetric.} Henceforth, we will {omit} the term $``$symmetric'' when referring to an object within $\Op_C(\SS)$ (or $\Coll_C(\SS)$).
\end{conv}

\begin{rem}\label{r:ulcat} An $\cS$-enriched category $\C$ can be identified with an $\cS$-enriched operad concentrated in arity $1$, i.e. defined by letting $\C(c;d) := \Map_\C(c,d)$. On the other hand, every $\cS$-enriched operad $\P$  has an \textbf{underlying category} $\P_1$ such that $\Map_{\P_1}(c,d) := \P(c;d)$. 
\end{rem}

{The following provides two typical examples of operads.}

\begin{define}  {(1) The \textbf{associative operad} $\Ass \in \Op_*(\cS)$ is defined by, for each $n\in\NN$, letting $\Ass(n) = \underset{\Sigma_n}{\bigsqcup}1_\cS$ with the symmetric-action being induced by the right regular representation of $\Sigma_n$ and with the composition induced by the concatenation of linear orders.}
	
	\smallskip
	
	{(2) The \textbf{commutative operad} $\Com\in \Op_*(\cS)$ is defined by letting $\Com(n) = 1_\cS$ for every $n\in\NN$, endowed with the trivial action by the symmetric groups, and such that the composition is induced by the canonical isomorphisms $1_\cS\otimes \cdots \otimes 1_\cS \x{\cong}{\lrar} 1_\cS$.}	
\end{define}

One can organize all the colored operads in $\SS$, for arbitrary sets of colors, into a single category as follows.

\begin{define} The \textbf{category of $\SS$-enriched operads}, denoted by $\Op(\SS)$, is the Grothendieck construction $$ \Op(\SS) := \int_{C\in \Sets} \Op_C(\SS) $$ where for each map $\alpha : C \lrar D$ of sets, the corresponding functor $\alpha^{*} : \Op_D(\SS) \lrar \Op_C(\SS)$ is defined by sending $\Q\in\Op_D(\SS)$ to $\alpha^{*}\Q$ with $$ \alpha^{*}\Q(c_1,\cdots,c_n;c) := \Q(\alpha(c_1),\cdots,\alpha(c_n);\alpha(c)).$$ 
	
	{More precisely, an object of $\Op(\SS)$ is the choice of a pair $(C,\P)$ with $C\in \Sets$ and $\P\in \Op_C(\SS)$. Moreover, a morphism $(C,\P) \lrar (D,\Q)$ consists of a map $\alpha : C \lrar D$ of sets and a map $f :  \P \lrar \alpha^{*}\Q$ in $\Op_C(\SS)$.}	
\end{define}

\subsection{Modules over an operad}\label{s:opmodules}

{Let $(\SS, \otimes, 1_\cS)$ be as in the previous subsection, and let $\P$ be a $C$-colored operad in $\cS$, for some fixed set of colors $C$}.

\smallskip

{Various types of operadic modules come naturally when considering $\P$ as a monoid in the monoidal category of $C$-collections, as expressed in the following.}

\begin{dfn} \label{d:operadicmodules}

	(i)	 {A \textbf{left $\mathcal{P}$-module} is a $C$-collection $M$ that is equipped with a left $\P$-action $\mathcal{P}\circ M\lrar M$ whose data consist of $\Sigma_*$-equivariant maps of the form
		\begin{gather*}
			\P(c_1,\cdots,c_n;c) \otimes M(d_{1,1},\cdots,d_{1,k_1};c_1) \otimes \cdots \otimes M(d_{n,1},\cdots,d_{n,k_n};c_n) \\
			\lrar M(d_{1,1},\cdots,d_{1,k_1},\cdots,d_{n,1},\cdots,d_{n,k_n};c).
		\end{gather*}
			These are required to satisfy the classical axioms of associativity and unitality for left modules.}
		
		\smallskip
		
		(ii) {Dually, a \textbf{right $\mathcal{P}$-module} is a $C$-collection $M$ equipped with a right $\mathcal{P}$-action $M \circ \P \lrar M$ whose data consist of $\Sigma_*$-equivariant maps of the form
			\begin{gather*}
				M(c_1,\cdots,c_n;c) \otimes \P(d_{1,1},\cdots,d_{1,k_1};c_1) \otimes \cdots \otimes \P(d_{n,1},\cdots,d_{n,k_n};c_n) \\
				\lrar M(d_{1,1},\cdots,d_{1,k_1},\cdots,d_{n,1},\cdots,d_{n,k_n};c)	
			\end{gather*}
			satisfying the classical axioms of associativity and unitality for right modules.}
		
		\smallskip
		
		(iii) {A $\P$\textbf{-bimodule} is a $C$-collection $M$ equipped with both a left and a right $\mathcal{P}$-module structure, which must satisfy the essential axiom of compatibility.} 
\end{dfn}  

We will write $\LMod(\P)$, $\RMod(\P)$ and $\BMod(\P)$ respectively for the categories of \textbf{left} $\P$-\textbf{modules}, \textbf{right} $\P$-\textbf{modules} and $\P$-\textbf{bimodules}.

\smallskip

The structure of right $\P$-modules is \textit{linear} in the sense that the forgetful functor $$\RMod(\P) \lrar \Coll_C(\SS)$$ preserves all colimits. However, an analogue in general does not hold for left $\P$-modules, and thus, does not hold for $\P$-bimodules as well. This promotes the development of \textit{infinitesimal left modules (bimodules)} over an operad (cf. \cite{Vallette}), essentially coming as the linearization of the former.

\smallskip

{We write $ - \circ_{(1)} - : \Coll_C(\cS) \times \Coll_C(\cS) \lrar \Coll_C(\cS)$ for the \textit{infinitesimal composite product} that is characterized as the \textit{right linearization} of the composite product $- \circ-$ (cf. \cite{Loday, Hoang}).} 

\begin{define}  \label{d:infmodules}

		(i) {An \textbf{infinitesimal left} $\P$-\textbf{module} is a $C$-collection $M$ endowed with a map $\P\circ_{(1)}M \lrar M$ whose data consist of $\Sigma_*$-equivariant maps of the form
			$$ \P(c_1,\cdots,c_n;c) \otimes M(d_1,\cdots,d_m;c_i) \lrar M(c_1,\cdots,c_{i-1},d_1,\cdots,d_m,c_{i+1},\cdots,c_n;c).$$
			These maps must satisfy the classical axioms of associativity and unitality for left modules.}
		
		\smallskip
		
		(ii) {An \textbf{infinitesimal} $\P$-\textbf{bimodule} is a $C$-collection $M$ equipped with both an infinitesimal left $\P$-module and a right $\P$-module structure. Moreover, these two structures are subject to the usual compatibility axiom. (See {\cite[$\S$2.2]{Hoang}}  for more details).}
\end{define}

We denote by $\IbMod(\P)$ (resp. $\ILMod(\P)$) the \textbf{category of infinitesimal $\P$-bimodules} (resp. \textbf{infinitesimal left $\P$-modules}). We will revisit these two notions in $\S$\ref{s:infbi}.

\begin{rem}\label{r:encode} It is noteworthy that each of the categories $\LMod(\P)$, $\BMod(\P)$, $\IbMod(\P)$ and $\RMod(\P)$ can be represented as the category of algebras over an enriched operad (or $\cS$-valued functors on an enriched category). For more details about these, {we refer the reader} to \cite{Hoang, Turchin, Julien}. Moreover, according to \cite{Guti}, there exists a discrete operad, {which we denote by}  $\textbf{O}_C$, that encodes the category of $C$-colored operads in $\cS$. 
\end{rem}

One of the main interests in the study of operads is the notion of algebras over an operad.

\begin{define} {A $\P$-\textbf{algebra} is a left $\P$-module concentrated in level $0$. More explicitly, a $\P$-algebra $A$ consists of a collection $\{A(c)\}_{c\in C}$ of objects in $\cS$ equipped, for each tuple $(c_1,\cdots,c_n;c)$, with an action map 
		$$ \P(c_1,\cdots,c_n;c)\otimes A(c_1)\otimes\cdots\otimes A(c_n) \lrar A(c).$$
		These maps must satisfy the essential axioms of associativity, unitality and equivariance. We denote by $\Alg_\P(\SS)$ the \textbf{category of $\P$-algebras}.}
\end{define}

\begin{example} {In line with the terminologies, the operad $\Ass$ (resp. $\Com$) encodes the category of associative monoids (resp. commutative monoids) in $\cS$.}
\end{example}

Next, the notion of modules over a $\P$-algebra $A$ essentially provides \textit{representations} for $A$.

\begin{dfn} An $A$\textbf{-module over} $\P$ is a collection  $M=\{M(c)\}_{c\in C}$ of objects in $\cS$ endowed, for each tuple $(c_1,\cdots,c_n;c)$ {and $k\in\{1,\cdots,n\}$}, with an action map of the form
	$$  \P(c_1,\cdots,c_n;c) \otimes \;  \bigotimes _{i \in \{1,\cdots,n\} \setminus \{k\}} \; A(c_i) \otimes M(c_k) \lrar M(c) .$$
	These maps are required to satisfy the essential axioms of associativity, unitality and equivariance. The \textbf{category of $A$-modules over} $\P$ will be denoted by $\Mod_\P^{A}$.
\end{dfn}

\begin{rem} As in \cite{BergerMoerdijk},  {$\P_A$  denotes} the \textbf{enveloping operad} of $A\in\Alg_\P(\cS)$. A key feature of this construction is the existence of a categorical isomorphism
	\begin{equation}\label{eq:AmodP}
		\Mod_\P^{A} \cong \Fun((\P_A)_1 , \cS)  
	\end{equation}
	between $A$-modules over $\P$ and $\cS$-valued enriched functors on the underlying category of $\P_A$ (cf. {\cite[Theorem 1.10]{BergerMoerdijk}}).
\end{rem}

In particular, when $\P$ is a single-colored operad then $(\P_A)_1$ forms an associative monoid in $\cS$. In this situation, we will write $\sU_\P(A) := (\P_A)_1$ and as usual, refer to it as the \textbf{universal enveloping algebra} of $A$ (cf. e.g., \cite{Loday, Fresse1}). The following is a fundamental result.   

\begin{prop}\label{p:uenv} Suppose that $\P$ is a single-colored operad. Then $\sU_\P(A)$ can be modeled by $\Free_A(1_{\cS})$ the free $A$-module over $\P$ generated by $1_{\cS}$.
	\begin{proof} As discussed above, the category $\Mod_\P^{A}$ is identified with the category of left modules over $\sU_\P(A)$. Note that $\sU_\P(A)$ is the free left module over itself that is generated by $1_{\cS}$. Thus we obtain an identification $\sU_\P(A)\cong\Free_A(1_{\cS})$ of $A$-modules over $\P$.
	\end{proof}
\end{prop}

\subsection{Semi-model categories}\label{s:semi}

{Let us now recall briefly some preliminaries relevant to semi-model categories. For more details, {we refer the reader} to \cite{Spitzweck, Fresse1}.}

\begin{define}{(\cite{Spitzweck}) The data of a \textbf{semi-model category} consist of a bicomplete category $\textbf{M}$ and together with three subcategories of \textit{weak equivalences}, \textit{fibrations} and \textit{cofibrations} such that the following conditions are satisfied.}
	
			\smallskip
	
	{(SM1) The class of weak equivalences satisfies the \textit{two-out-of-three} axiom.}
	
			\smallskip
	
	{(SM2) The three classes of morphisms are all closed under retracts.}
	
			\smallskip
	
	{(SM3) (i) The cofibrations have the left lifting property with respect to \textit{trivial fibrations}; and (ii) the \textit{trivial cofibrations} with \textit{cofibrant} domains have the left lifting property with respect to fibrations. Here, a morphism in $\textbf{M}$ is called a trivial (co)fibration if it is both a weak equivalence and a (co)fibration; and an object $X \in \textbf{M}$ is called cofibrant if the unique morphism $\emptyset_{\textbf{M}} \lrar X$ is a cofibration, where $\emptyset_{\textbf{M}}$ denotes the initial object of $\textbf{M}$.}
	
			\smallskip
	
	{(SM4) (i) Every morphism in $\textbf{M}$ can be functorially factored into a cofibration followed by a
		trivial fibration, and (ii) every morphism in $\textbf{M}$ whose domain is cofibrant can be functorially factored into a trivial cofibration followed by a fibration.}
	
			\smallskip
	
	{(SM5) The two classes of fibrations and trivial fibrations are closed under pullbacks.}
	
			\smallskip
	
	{(SM6) The initial object $\emptyset_{\textbf{M}}$ is cofibrant.}
\end{define}

\begin{define}{(\cite{Spitzweck}) A semi-model category $\textbf{M}$ is said to be \textbf{cofibrantly generated} if there are sets of morphisms $I$ and $J$ such that the following conditions are fulfilled:}
	
	\smallskip
	
	{(i) The class of morphisms which have the right lifting property with respect to $J$ (resp. $I$) coincides with that of fibrations (resp. trivial fibrations).}
	
	\smallskip
	
	{(ii) The domains of maps in $I$ are small relative to $I$-\textit{cell}, and the domains of maps in $J$ are
		small relative to maps in $J$-\textit{cell} whose domains are cofibrant. Here, $I$-cell (resp. $J$-cell) refers to the class of maps which are transfinite compositions of pushouts of maps from $I$ (resp. $J$).}
	
	\smallskip
	
	{We will refer to $I$ (resp. $J$) as the set of \textit{generating cofibrations} (resp. \textit{trivial cofibrations}).}
	
	\smallskip
	
	{Moreover, a semi-model category $\textbf{M}$ is \textbf{combinatorial} if it is cofibrantly generated, and in addition, the underlying category $\textbf{M}$ is locally presentable.}
\end{define}

\begin{rem} {As in (full) model categories, the \textit{homotopy category} $\Ho(\textbf{M})$ of a semi-model category $\textbf{M}$ can be modeled by $\pi\textbf{M}_{cf}$ the category whose objects are \textit{bifibrant} (i.e. both fibrant and cofibrant) objects and whose morphisms are given by the homotopy classes of maps in $\textbf{M}$. Moreover, one may also define \textit{function complexes} between objects using \textit{(co)simplicial frames} of objects in $\textbf{M}$. As usual, these serve as models for the \textit{derived mapping space} $\Map^{\h}_{\textbf{M}}(-,-)$ in $\textbf{M}$. For more details, {we refer the reader} to \cite{Spitzweck} once again.}
\end{rem}

\begin{define} An adjunction $\L : \adjunction*{}{\textbf{M}}{\textbf{N}}{} : \R$ between semi-model categories is called a \textbf{Quillen adjunction} if the right adjoint $\R$ preserves fibrations and trivial fibrations.
\end{define}

For a Quillen adjunction $\L \dashv \R$ between semi-model categories, one can show that the left adjoint $\L$ (resp. the right adjoint $\R$) preserves weak equivalences between cofibrant (resp. bifibrant) objects, due to a version of \textit{K. Brown's lemma}. In light of this, the adjunction $\L \dashv \R$ induces an adjunction $\mathbb{L}\L : \adjunction*{}{\Ho(\textbf{M})}{\Ho(\textbf{N})}{} : \mathbb{R}\R$ between homotopy categories by taking the corresponding \textit{left and right derived functors}. 

\begin{define} A Quillen adjunction $\L : \adjunction*{}{\textbf{M}}{\textbf{N}}{} : \R$ between semi-model categories is a \textbf{Quillen equivalence} if the induced adjunction $\mathbb{L}\L \dashv \mathbb{R}\R$ is an equivalence between homotopy categories.
\end{define}

The reader may refer to \cite{Fresse1} for various equivalent characterizations of a Quillen equivalence between semi-model categories, which are quite analogous to those in the context of (full) model categories.

\subsection{Operadic model structures}\label{s:optransfermod}

{We assume  that} $(\SS, \otimes, 1_\cS)$ is a \textbf{symmetric monoidal model category} (cf. \cite{Hovey}). Let $\P$ be a $C$-colored operad in $\cS$, {where $C$ is a fixed set of colors}. We will usually require the existence of a \textit{transferred model (or semi-model) structure} on the category $\Alg_\P(\cS)$. This is based on the  result {below. By definition, we say that $\P$ is $\Sigma$-\textbf{cofibrant} if it is cofibrant as an object in $\Coll_C(\cS)$ with respect to the  model structure transferred from $\cS$.}

\begin{thm}\label{t:semi}(Spitzweck \cite{Spitzweck}, Fresse \cite{Fresse1}) {Suppose  that} the model structure on $\cS$ is cofibrantly generated, that the monoidal unit $1_\cS$ is cofibrant, and that $\P$ is $\Sigma$-cofibrant. Then the category $\Alg_\P(\cS)$ carries a cofibrantly generated semi-model structure transferred from the model structure on $\cS^{\times C}$. {Namely, a map $A \lrar B$ in $\Alg_\P(\cS)$ is a weak equivalence (resp. fibration) if and only if for every $c\in C$ the  map $A(c) \lrar B(c)$ is a weak equivalence (resp. fibration) in $\cS$.} 
\end{thm} 

\begin{rem} In particular, with the same assumptions on $\cS$, the transferred semi-model structure on the category $\Op_C(\SS)$ does exist. Additionally, with the same assumptions on $\cS$ and $\P$, {the categories} $\LMod(\P)$ and $\BMod(\P)$ admit a transferred semi-model structure as well, while $\RMod(\P)$ and $\IbMod(\P)$ inherit a transferred (full) model structure. These observations follow from the above theorem and Remark \ref{r:encode}.
\end{rem}

We now recall the {Dwyer-Kan model structure on} $\cS$-enriched operads, according to the work of Caviglia \cite{Caviglia}. For $\C \in \Cat(\cS)$ an $\cS$-\textbf{enriched category}, the \textbf{homotopy category} of $\C$, denoted $\Ho(\C)$, is the ordinary category whose objects are the same as those of $\C$ and such that, for $x,y\in \Ob(\C)$, we have  $$\Hom_{\Ho(\C)}(x,y) := \Hom_{\Ho(\cS)}(1_\cS,\Map_\C(x,y)) .$$
By convention, the \textbf{homotopy category of an operad} $\P$ is $\Ho(\P) := \Ho(\P_1)$.  
Moreover, a map $f : \P \rar \Q$ in $\Op(\cS)$ is called a \textit{levelwise weak equivalence} (resp. \textit{fibration, trivial fibration}, etc) if for every sequence $(c_1,\cdots,c_n;c)$ of colors of $\P$, the induced map $$\P(c_1,\cdots,c_n;c) \lrar \Q(f(c_1),\cdots,f(c_n) ; f(c))$$ is a weak equivalence (resp. fibration, trivial fibration, etc) in $\cS$.

\begin{dfn}\label{d:DKop} A map $f : \P \lrar \Q$ in $\Op(\cS)$ is called a \textbf{Dwyer-Kan equivalence} if it is a levelwise weak equivalence and such that the induced functor $$\Ho(f):\Ho(\P)\lrar\Ho(\Q)$$ between homotopy categories is essentially surjective.
\end{dfn}

\begin{rem}\label{d:DKope} {Under mild assumptions on $\cS$ (in particular, combinatoriality), there is a combinatorial \textbf{Dwyer-Kan model structure} on $\Op(\cS)$, in which the weak
		equivalences are the Dwyer-Kan equivalences and the trivial fibrations are the
		levelwise trivial fibrations that are surjective on colors. For more details, we refer the reader to   \cite[$\S$8.6]{Caviglia}.}
\end{rem}

\begin{examples}\label{ex:base} Here are some base categories that we are concerned with in this paper. Let $\textbf{k}$ be a (unital) commutative ring.
	
	\smallskip
	
	(i) Firstly, we let $\Set_\Delta$ denote the category of \textit{simplicial sets}. This category comes equipped with the Cartesian monoidal structure and with the standard (Kan-Quillen) model structure. 
	
	\smallskip
	
	(ii) Another instance is the category of \textit{simplicial} $\textbf{k}$-\textit{modules}, denoted $\sMod(\textbf{k})$, and equipped with the (degreewise) tensor product over $\textbf{k}$, and with the standard model structure transferred from that on $\Set_\Delta$. 
	
	\smallskip
	
	(iii) We will denote by $\C(\textbf{k})$ (resp. $\C_{\geqslant0}(\textbf{k})$) the category of \textit{dg $\textbf{k}$-modules} (resp. \textit{connective dg $\textbf{k}$-modules}), endowed with the usual tensor product, and with the projective model structure. 
\end{examples}

\begin{rem} When $\cS = \Set_\Delta$ or $\sMod(\textbf{k})$, then for every operad $\P$ the transferred homotopy theory on $\Alg_\P(\cS)$ exists as a (full) model category. {Moreover, the Dwyer-Kan model structure on $\Op(\cS)$ does exist as well.}
\end{rem}

\begin{rem}\label{r:Ck} For the case $\cS = \C(\textbf{k})$ or $\C_{\geqslant0}(\textbf{k})$, due to Theorem \ref{t:semi} the transferred homotopy theory on $\Alg_\P(\cS)$ exists as a semi-model category,  provided that $\P$ is $\Sigma$-cofibrant. Additionally, it can be shown that {the same classes of weak equivalences and generating (trivial) cofibrations as in {\cite[Theorem 8.6]{Caviglia}}  define a cofibrantly generated semi-model structure on $\Op(\cS)$. We will refer to this as the \textit{Dwyer-Kan semi-model structure} on $\Op(\cS)$}. (This is in fact made up together by the Dwyer-Kan model structure on $\Cat(\cS)$ and by the transferred semi-model structure on the categories $\Op_C(\cS)$ for $C\in\Sets$). In particular, in this situation $\cS$ is not \textit{sufficient} in the sense of {\cite[$\S$3]{Hoang}}. In spite of this, the results needed from {that paper} remain valid if we bring them in the framework of semi-model categories.
\end{rem}

\subsection{Spectral cotangent complex and  Quillen cohomology}\label{s:tangentcategory}

We recall briefly some fundamental notions relevant to the \textit{spectral cotangent complex} and \textit{Quillen cohomology}, according to the works of Harpaz-Nuiten-Prasma \cite{YonatanBundle, YonatanCotangent}. Nevertheless, our setting is a bit more flexible, so that we may define \textit{tangent categories} of (semi) model categories of interest without requiring the left properness. This enables us to cover a wider range of examples. {We also refer the reader to \cite[$\S$2.5]{Hoang}} for more details.

\smallskip

By definition, a {semi-model} category $\textbf{M}$ is \textbf{weakly pointed} if it contains a \textbf{weak zero object} $0$, i.e. $0$ is both homotopy initial and terminal. When $\textbf{M}$ is weakly pointed, we may define the \textit{suspension} and \textit{desuspension} of an object $X\in\Ob(\textbf{M})$ respectively as 
$$ \Sigma(X) := 0 \, \x{\h}{\underset{X}{\bigsqcup}} \, 0  \; \; \; \text{and} \; \; \; \Omega(X) := 0 \, \x{\h}{\underset{X}{\times}} \, 0. $$
Moreover, when passing to the homotopy category, these constructions determine an adjunction $$\Sigma : \adjunction*{}{\Ho(\textbf{M})}{\Ho(\textbf{M})}{}  : \Omega.$$

\begin{dfn}\label{d:stable} A {semi-model} category $\textbf{M}$ is said to be \textbf{stable} if the following equivalent conditions hold:
	
	\smallskip
	
	(1) The \textbf{underlying} $\infty$-\textbf{category} $\textbf{M}_\infty$ is stable in the sense of \cite{Lurieha}.
	
	\smallskip
	
	(2) $\textbf{M}$ is weakly pointed and such that a square in $\textbf{M}$ is homotopy coCartesian if and only if it is homotopy Cartesian.
	
	\smallskip
	
	(3) $\textbf{M}$ is weakly pointed and such that the self-adjunction $\Sigma \dashv \Omega$ on $\Ho(\textbf{M})$ is an adjoint equivalence.
\end{dfn}

We now recall the procedure of taking \textit{stabilizations}. Let $\textbf{M}$ be a \textbf{weakly pointed {semi-model} category} and let $X$ be an ($\mathbb{N}\times\mathbb{N} $)-diagram in $\textbf{M}$. A square of $X$ that takes the form
$$ \xymatrix{
	X_{n,n} \ar[r]\ar[d] & X_{n,n+1} \ar[d] \\
	X_{n+1,n} \ar[r] & X_{n+1,n+1} \\
}$$
will be called a \textit{diagonal square}.
\begin{dfn}\label{dfnspectrumobj} An  ($\mathbb{N}\times\mathbb{N} $)-diagram in $\textbf{M}$ is called
	
	\smallskip
	
	(1) a \textbf{prespectrum} if all its off-diagonal entries are weak zero objects in $\textbf{M}$,
	
	\smallskip
	
	(2) an \textbf{$\Omega$-spectrum} if it is a prespectrum and all its diagonal squares are homotopy Cartesian,
	
	\smallskip
	
	(3) a \textbf{suspension spectrum} if it is a prespectrum and all its diagonal squares are homotopy coCartesian.
\end{dfn}

\begin{dfn} A map $f : X \lrar Y$ in $\textbf{M}^{\mathbb{N}\times\mathbb{N}}$ is called a \textbf{stable equivalence} if for every $\Omega$-spectrum $Z$ the induced map 
	$$ \Map^{\der}_{\textbf{M}^{\mathbb{N}\times\mathbb{N}}_{\proj}}(Y,Z) \lrar  \Map^{\der}_{\textbf{M}^{\mathbb{N}\times\mathbb{N}}_{\proj}}(X,Z) $$
	is a weak equivalence, where $\textbf{M}^{\mathbb{N}\times\mathbb{N}}_{\proj}$ signifies the projective {semi-model} category of ($\mathbb{N}\times\mathbb{N} $)-diagrams in $\textbf{M}$.
\end{dfn}

We will assume further that $\textbf{M}$ is a \textit{weakly pointed combinatorial {semi-model} category} such that the domains of generating cofibrations are cofibrant.

\begin{dfn}\label{d:stabilization} The \textbf{stabilization} of $\textbf{M}$, denoted $\Sp(\textbf{M})$, is the left Bousfield localization of $\textbf{M}^{\mathbb{N}\times\mathbb{N}}_{\proj}$ with $\Omega$-spectra as the local objects. More precisely, $\Sp(\textbf{M})$ is a cofibrantly generated semi-model category such that
	
	\smallskip
	
	$\bullet$ weak equivalences are the stable equivalences, and 
	
	\smallskip
	
	$\bullet$ (generating) cofibrations are the same as those of $\textbf{M}^{\mathbb{N}\times\mathbb{N}}_{\proj}$.
	
	\smallskip
	
	\noindent In particular, fibrant objects of $\Sp(\textbf{M})$ are precisely the levelwise fibrant $\Omega$-spectra.
\end{dfn}

\begin{rem} The above definition is {valid due to} the main result of \cite{White}, and together with the fact that the $\Omega$-spectra in $\textbf{M}$ can be characterized as local objects against a certain set of maps (cf. \cite{YonatanBundle}, Lemma 2.1.6). 
\end{rem}
By construction, there is a Quillen adjunction $\Sigma^{\infty} : \adjunction*{}{\textbf{M}}{\Sp(\textbf{M})}{} : \Omega^{\infty}$ where $\Omega^{\infty}(X)=X_{0,0}$ and $\Sigma^{\infty}({A})$ is the constant diagram with value ${A}$.

\begin{notns} {For an object $A\in\Ob(\textbf{M})$, we will write $\textbf{M}_{/A}$ (resp. $\textbf{M}_{A/}$) for the \textit{category of objects over} (resp. \textit{under}) $A$. Moreover, we will denote by $\textbf{M}_{A//A} := (\textbf{M}_{/A})_{\Id_A/}$ and refer to it the \textit{category of objects over and under $A$}. For more information, an object of $\textbf{M}_{A//A}$ is a sequence of maps in $\textbf{M}$ of the form $A \x{g}{\lrar} B \x{f}{\lrar} A $ such that $f\circ g = \Id_A$.}  
\end{notns}

\begin{dfn}\label{d:tangentcat} For an object $A\in\Ob(\textbf{M})$, the \textbf{tangent category} to $\textbf{M}$ at $A$ is
	$$ \mathcal{T}_A\textbf{M}:=\Sp(\textbf{M}_{A//A})  $$
	the stabilization of the {semi-model} category of objects over and under $A$.	
\end{dfn}

There is a Quillen adjunction  $ \adjunction{\Sigma^{\infty}_+}{\textbf{M}_{/A}}{\mathcal{T}_A\textbf{M}}{\Omega^{\infty}_+} $ defined as the composed adjunction
$$ \adjunction*{A\sqcup (-)}{\textbf{M}_{/A}}{\textbf{M}_{A//A}}{\forget} \adjunction*{\Sigma^{\infty}}{}{\mathcal{T}_A\textbf{M}}{\Omega^{\infty}} .$$
For more details, the left adjoint $\Sigma^{\infty}_+$ sends each $B \in \textbf{M}_{/A}$ to $$ \Sigma^{\infty}_+(B) = \Sigma^{\infty}(A\longrightarrow A\sqcup B \longrightarrow A),$$ while the right adjoint $\Omega^{\infty}_+$ sends $X \in \mathcal{T}_A\textbf{M}$ to $\Omega^{\infty}_+(X) = (X_{0,0} \lrar A)$.

\smallskip

Here is the central concept of the paper.
\begin{dfn} For an object $A\in\Ob(\textbf{M})$, the \textbf{cotangent complex} of $A$ is 
	$$\rL_A:=\mathbb{L}\Sigma^{\infty}_+(A) \in  \mathcal{T}_A\textbf{M},$$
	i.e. the derived suspension spectrum of $A$. 
\end{dfn}

\begin{rem}\label{r:sus} Accordingly, the cotangent complex of $A$ can be modeled by the constant spectrum $\Sigma^{\infty}(A\sqcup A^{\cof})$ where $A^{\cof}$ signifies a cofibrant resolution of $A\in\textbf{M}$. Furthermore, $\rL_A = \Sigma^{\infty}(A\sqcup A^{\cof})$ admits a suspension spectrum replacement determined by fixing $A\sqcup A^{\cof}$ as its value at the bidegree $(0,0)$. In particular, the value at the bidegree $(n,n)$ is given by $\Sigma^n(A\sqcup A^{\cof})$ the $n$-suspension of $A\sqcup A^{\cof}\in \textbf{M}_{A//A}$ (see {
		\cite[Corollary 2.3.3]{YonatanBundle}}). In practice, it is usually convenient to exhibit that suspension spectrum as a model for $\rL_A$. 
\end{rem}

\begin{dfn}\label{d:relcotangent} For a map $f : A \lrar B$ in $\textbf{M}$, the \textbf{relative cotangent complex} of $f$ is 
	$$\rL_{B/A}:=\hocofib \, [\,\mathbb{L}\Sigma^{\infty}_+(f) \lrar \rL_B\,]$$
	the homotopy cofiber of the map $\mathbb{L}\Sigma^{\infty}_+(f) \lrar \rL_B$ in $\mathcal{T}_B\textbf{M}$.
\end{dfn}

\begin{rem} When $A$ is the initial object then $\rL_{B/A}$ is weakly equivalent to $\rL_B$. Thus, on the one hand, the relative cotangent complex can be viewed as a natural extension of the standard notion of cotangent complex. On other hand, for any map $f : A \lrar B$, we can identify $\rL_{B/A}$ with the (absolute) cotangent complex of $f$ considered as an object of $\textbf{M}_{A/}$ (see {\cite[Proposition 2.2.10]{YonatanCotangent}}). 
\end{rem}

\begin{dfn}\label{d:Qcohom} Let $X$ be a fibrant object of $\textbf{M}$. The \textbf{(spectral) Quillen cohomology} of $X$ with coefficients in a given object $M\in\mathcal{T}_X\textbf{M}$ is the space
	$$ \HHQ^\star(X ; M) :=  \Map^{\der}_{\mathcal{T}_X\textbf{M}} (\rL_X , M).$$
	Moreover, for each $n\in\mathbb{Z}$, the \textit{$n$-th Quillen cohomology group of $X$ with coefficients in $M$} is given by
	$$ \HHQ^n(X  ; M) := \pi_0 \Map^{\der}_{\mathcal{T}_X\textbf{M}} (\rL_X , M[n])$$
	where $M[n] : = \Sigma^{n}M$ the $n$-suspension of $M$ in $\mathcal{T}_X\textbf{M}$.
\end{dfn}

\begin{rem}\label{r:trivclass} We will regard $\HHQ^\star(X ; M)$ as a pointed space with the base point being the zero map $\rL_X \x{0}{\lrar} M$. Note also that this base point represents the zero element of the group $\HHQ^0(X  ; M)$.
\end{rem}

\subsection{More conventions and notations}\label{s:not}

Let $\cS$ be a suitable symmetric monoidal model category, $\P$ a $C$-colored operad in $\cS$ {with $C$  some fixed set of colors}, and let $A$ be a $\P$-algebra. To proceed, we will use the following notations and conventions.

\smallskip

(i) We will say that the pair $(\cS, \P)$ is \textbf{sufficient} (resp. \textbf{stably sufficient}) if the base category $\cS$ is \textit{sufficient} (resp. \textit{stably sufficient}) in the sense of {\cite[$\S$3]{Hoang}}, and in addition, $\P$ is $\Sigma$-cofibrant.

\smallskip

(ii) Moreover, the pair $(\cS, \P)$ is said to be \textbf{abundant} (resp. \textbf{stably abundant}) if the base category $\cS$ is \textit{abundant} (resp. \textit{stably abundant}) in the sense of {\cite[$\S$3]{Hoang}} and such that $\P$ is cofibrant, and is \textbf{good} in the sense of {\cite[$\S$5.1]{Hoang}}.

\smallskip

(iii) We will say that the triple $(\cS,\P, A)$ is \textbf{sufficient} (resp. \textbf{stably sufficient}) if the pair $(\cS, \P)$ is \textit{sufficient} (resp. \textit{stably sufficient}) in the above sense, and in addition, $A\in \Alg_\P(\cS)$ is cofibrant. 

\smallskip

(iv) Furthermore, $(\cS,\P, A)$ is \textbf{abundant} (resp. \textbf{stably abundant}) if the pair $(\cS, \P)$ is \textit{abudant } (resp. \textit{stably abundant}) and $A\in \Alg_\P(\cS)$ is cofibrant.

\smallskip

The hypotheses presented above allow us to inherit some main results of \cite{Yonatan, Hoang}. Therefore, readers may safely disregard the details of these hypotheses.

\smallskip

(v)  {As mentioned in $\S$\ref{s:opalg}}, we will write $\P^{\sB}$ (resp. $\P^{\Inf}$ and $\P^{\sr}$) to denote $\P$ in its role as an object of $\BMod(\P)$ (resp. $\IbMod(\P)$ and $\RMod(\P)$), and write {$A^{\me}$ for $A$ as an object in $\Mod_\P^A$.}

\smallskip

(vi) We will refer to different cotangent complexes of $\P$ using the following notations:

\smallskip

\; $\bullet$ $\rL_\P \in \T_\P\Op(\cS)$: the cotangent complex of $\P$ considered as an object of $\Op(\cS)$.

\smallskip

\; $\bullet$ $\rL^{\red}_\P \in \T_\P\Op_C(\cS)$: the cotangent complex of $\P$ considered as an object of $\Op_C(\cS)$.

\smallskip

\; $\bullet$ $\rL_{\P^{\sB}} \in \T_{\P^{\sB}}\BMod(\P)$: the cotangent complex of $\P$ when regarded as a bimodule over itself.

\smallskip

\; $\bullet$ $\rL_{\P^{\Inf}} \in \T_{\P^{\Inf}}\IbMod(\P)$: the cotangent complex of $\P$ as an infinitesimal bimodule over itself.

\smallskip

(vii) To denote different cotangent complexes of $A$, we will use the following:    

\smallskip

\; $\bullet$ $\rL_A \in \T_A\Alg_\P(\cS)$: the cotangent complex of $A$ considered as an algebra over $\P$, and

\smallskip

\; $\bullet$ $\rL_{A^{\me}} \in \T_{A^{\me}} \Mod_\P^A$: the cotangent complex of $A$ considered as an object of $\Mod_\P^A$.

\smallskip

(viii) {We let $\HHQ^\star(\P; -) $ and $ \HHQ^\star_{\red}(\P; -)$ denote the Quillen cohomology of $\P$ as an object of $\Op(\cS)$, and $\Op_C(\cS)$, respectively. The former will be called the (\textit{proper}) \textit{Quillen cohomology of} $\P$ (classified by $\rL_\P \in \T_\P\Op(\cS)$), and the latter \textit{reduced} (classified by $\rL^{\red}_\P \in \T_\P\Op_C(\cS)$).}

\smallskip

{Finally, $\HHQ^\star(A; -)$ denotes the Quillen cohomology of $A$ as a $\P$-algebra.}

\smallskip

{(ix)} We will write $$\Free_\P^{\Inf} : \Coll_C(\cS) \lrar \IbMod(\P) \;\;\; (\text{resp.} \, \Free_\P^{i\ell} : \Coll_C(\cS) \lrar \ILMod(\P))$$ 
for the functor that assigns to each $C$-collection $M$ the \textbf{free infinitesimal $\P$-bimodule} (resp. \textbf{free infinitesimal left $\P$-module}) generated by $M$.

\smallskip

{(x)} We will write $\Free_A : \cS^{\times C} \lrar \Mod_A^\P$ {for} the functor that assigns to each object $X\in \cS^{\times C}$ the free $A$-module over $\P$ generated by $X$.

\smallskip

{(xi)} We will denote by $\mathcal{E}_C$ the $C$-collection that is $\emptyset_\cS$ on all levels except that $\mathcal{E}_C(c) = 1_\cS$ for every $c\in C$. When $C$ is a singleton, {this is denoted by} $\mathcal{E}_*$.

\section{Operadic tangent categories and cotangent complex of enriched operads}\label{s:rv}

We recall from  \cite{Yonatan, Hoang} some fundamental results concerning various operadic tangent categories and cotangent complex of operads. Let $\cS$ be a symmetric monoidal model category, $\P$ a $C$-colored operad in $\cS$ {with $C$  some fixed set of colors}, and let $A$ be a $\P$-algebra. 

\smallskip

The following theorem plays a key role in the investigation of operadic tangent categories.

\begin{thm}\label{t:tanalg}(Harpaz-Nuiten-Prasma \cite{Yonatan}) {Suppose that the triple $(\cS,\P,A)$ is \textit{sufficient}.} There is a Quillen equivalence
	\begin{equation}\label{eq:Quillenalg}
		\adjunction*{\simeq}{\T_{A^{\me}} \Mod_\P^A}{\T_A\Alg_\P(\cS)}{}
	\end{equation}
	induced by the adjunction $\adjunction*{}{(\Mod_\P^A)_{A^{\me}/}}{(\Alg_\P(\cS))_{A/}}{}$ of induction-restriction functors. Moreover, when $(\cS,\P,A)$ is \textit{stably sufficient} then there is a chain of Quillen equivalences
	\begin{equation}\label{eq:Quillenalg1}
		\adjunction*{\simeq}{\Mod_\P^A}{}{}\adjunction*{\simeq}{\T_{A^{\me}} \Mod_\P^A}{\T_A\Alg_\P(\cS)}{}
	\end{equation}
	where the first adjunction is given by the composed Quillen equivalence
	\begin{equation}\label{eq:seqQuillenstab1}
		\adjunction*{(-) \sqcup A^{\me}}{\Mod_\P^A}{(\Mod_\P^A)_{A^{\me}//A^{\me}}}{\ker} \adjunction*{\Sigma^{\infty}}{}{\T_{A^{\me}} \Mod_\P^A}{\Omega^{\infty}}
	\end{equation}
	(see $\S$\ref{s:not} for notations).
\end{thm}

A version of the above theorem for various tangent categories at $\P$ is as follows.

\begin{thm}\label{t:tanop}(\cite{Hoang}) Assume that the pair $(\cS, \P)$ is \textit{sufficient}. There is a chain of Quillen equivalences
	\begin{equation}\label{eq:seqQuillen}
		\adjunction*{\simeq}{\T_{\P^{\Inf}}\IbMod(\P)}{\T_{\P^{\sB}}\BMod(\P)}{} \adjunction*{\simeq}{}{\T_\P\Op_C(\mathcal{S})}{} \adjunction*{\simeq}{}{\T_\P\Op(\mathcal{S})}{}
	\end{equation}
	that is induced by the adjunctions of induction-restriction functors
	$$ \adjunction*{}{\IbMod(\P)_{\P^{\Inf}/}}{\BMod(\P)_{\P^{\sB}/}}{} \adjunction*{}{}{\Op_C(\mathcal{S})_{\P/}}{} \adjunction*{}{}{\Op(\mathcal{S})_{\P/}}{}.$$ Furthermore, when $(\cS, \P)$ is \textit{stably sufficient}, we have a prolonged chain of Quillen equivalences:
	\begin{equation}\label{eq:seqQuillenstab}
		\adjunction*{\simeq}{\IbMod(\P)}{}{}\adjunction*{\simeq}{\T_{\P^{\Inf}}\IbMod(\P)}{\T_{\P^{\sB}}\BMod(\P)}{} \adjunction*{\simeq}{}{\T_\P\Op_C(\mathcal{S})}{} \adjunction*{\simeq}{}{\T_\P\Op(\mathcal{S})}{}
	\end{equation}
	where the first adjunction is given by the composed Quillen equivalence
	\begin{equation}\label{eq:seqQuillenstab2}
		\adjunction*{(-) \sqcup \P^{\Inf}}{\IbMod(\P)}{\IbMod(\P)_{\P^{\Inf}//\P^{\Inf}}}{\ker} \adjunction*{\Sigma^{\infty}}{}{\T_{\P^{\Inf}}\IbMod(\P)}{\Omega^{\infty}}
	\end{equation}
	(see $\S$\ref{s:not} for notations).
\end{thm}

In what follows, we recall from \cite{Hoang} a description of the cotangent complex of $\P$. {The following is the key proposition.}

\begin{prop}\label{p:key}{(\cite[$\S$5.2]{Hoang})} Suppose that the pair $(\cS, \P)$ is \textit{abundant}. Under the right Quillen equivalence $\T_\P\Op(\mathcal{S}) \lrarsimeq \T_{\P^{\sB}}\BMod(\P)$,
	the cotangent complex $\rL_\P \in \T_\P\Op(\cS)$ is identified with $\rL_{\P^{\sB}}[-1] \in \T_{\P^{\sB}}\BMod(\P)$. 
\end{prop}

\begin{notns} \label{no:sphere}
		(i) For each $n \in \NN$, we denote by $\sS^{n}:=\Sigma^{n}(1_\mathcal{S}\sqcup 1_\mathcal{S})\in\cS$ with the suspension $\Sigma(-)$ taken in $\mathcal{S}_{1_\mathcal{S}//1_\mathcal{S}}$, and refer to $\sS^{n}$ as the \textbf{pointed} $n$\textbf{-sphere} in $\cS$.
		
		\smallskip
		
		(ii) We will write $\sS_C^{n}$ to denote the $C$-collection that is $\emptyset_\cS$ on all levels except that $\sS_C^{n}(c;c)=\sS^{n}$ for every $c\in C$.
		
		\smallskip
		
		(iii) {Furthermore, we let $\P\circ\sS_C^{\bullet}\in\T_{\P^{\sB}}\BMod(\P)$ denote the suspension spectrum with $(\P\circ\sS_C^{\bullet})_{n,n} = \P\circ\sS_C^{n}$ for every $n \in \NN$. (This agrees with the object $\widetilde{\rL}_\P$  {from \cite[$\S$5.2]{Hoang}}). In particular, we have on each level that
			$$ (\P\circ\sS_C^{n})(c_1,\cdots,c_m;c) = \P(c_1,\cdots,c_m;c)\otimes (\sS^{n})^{\otimes m}.$$}	
\end{notns}

\begin{rem} {For each $n \in \NN$, we have a canonical weak equivalence $\P\circ\sS_C^{n} \simeq \Sigma^n(\P^{\sB} \sqcup \P^{\sB})$ in $\BMod(\P)_{\P^{\sB}//\P^{\sB}}$. Consequently, $\P\circ\sS_C^{\bullet}$ is a suspension-spectrum model for $\rL_{\P^{\sB}}\in\T_{\P^{\sB}}\BMod(\P)$.}
\end{rem}

In summary, the cotangent complex $\rL_\P \in \T_\P\Op(\cS)$ is described as follows.

\begin{thm}(\cite{Hoang}) \label{t:maincot} Suppose that the pair $(\cS, \P)$ is \textit{abundant}. Then under the right Quillen equivalence
	$$ \T_\P\Op(\cS) \lrarsimeq \T_{\P^{\Inf}}\IbMod(\P),$$
	the cotangent complex $\rL_\P$ is identified with $(\P\circ\sS_C^{\bullet})[-1]\in\T_{\P^{\Inf}}\IbMod(\P)$. {(Here, we use the same notation for the underlying prespectrum in $\T_{\P^{\Inf}}\IbMod(\P)$ of $\P\circ\sS_C^{\bullet}\in\T_{\P^{\sB}}\BMod(\P)$). Moreover, when $(\cS, \P)$ is \textit{stably abundant} then  under the right Quillen equivalence
		$$ \T_\P\Op(\cS) \lrarsimeq \IbMod(\P),$$
		$\rL_\P \in \T_\P\Op(\cS)$ is identified with $\ovl{\rL}_\P[-1]\in\IbMod(\P)$ where $\ovl{\rL}_\P \in \IbMod(\P)$ is given on each level by 
		\begin{equation}\label{eq:barL}
			\ovl{\rL}_\P(c_1,\cdots,c_m;c) \simeq \P(c_1,\cdots,c_m;c)\otimes \hocolim_n\Omega^{n} [\, (\sS^{n})^{\otimes m}  \times_{1_\mathcal{S}}^{\h}  0\,]
		\end{equation}
		with the desuspension $\Omega(-)$ taken in $\cS$.}
\end{thm}

\section{More about infinitesimal bimodules over an operad}\label{s:infbi}

{Let $(\SS, \otimes, 1_\cS)$ and $\P \in \Op_C(\cS)$ be as in $\S$\ref{s:opmodules}.} First recall from  {\cite[$\S$2.2]{Hoang}} that, for each $C$-collection $M$, we have
\begin{equation}\label{eq:freeinf}
	\Free_\P^{\Inf}(M) \cong \mathcal{P}\circ_{(1)}( M \circ \mathcal{P}) \;\;\; \text{and} \;\;\; \Free_\P^{i\ell}(M) \cong \mathcal{P}\circ_{(1)}M.
\end{equation}
Taking $M = \mathcal{E}_C$, we have
$$ \Free_\P^{i\ell}(\mathcal{E}_C) \cong \mathcal{P}\circ_{(1)}\mathcal{E}_C \cong \mathcal{P}\circ_{(1)}( \mathcal{E}_C \circ \mathcal{P}) \cong \Free_\P^{\Inf}(\mathcal{E}_C)$$
(see $\S$\ref{s:not} for notations). In particular, $\mathcal{P}\circ_{(1)}\mathcal{E}_*$ agrees with the $``$\textit{shifted object}'' $\P[1]$ of {\cite[$\S$17.3]{Fresse1}}, which is given on each level by $\P[1](n)= \mathcal{P}(n+1)$. Thus we may write
\begin{equation}\label{eq:freeinf1}
	\Free_\P^{i\ell}(\mathcal{E}_*) \cong  \P[1] \cong \Free_\P^{\Inf}(\mathcal{E}_*).
\end{equation}

\begin{example}\label{r:comfree} A basic fact is that $\P$ is free (generated by $\I_C$) when regarded as a right module over itself. Nonetheless, $\P$ is in general not a free object when living in $\ILMod(\P)$ (or $\IbMod(\P)$){, except in the case of the commutative operad $\Com$. Indeed, due to \eqref{eq:freeinf1}} the operad $\Com$ is free (generated by $\mathcal{E}_*$) when considered as an object of $\IbMod(\Com)$ (or $\ILMod(\Com)$).
\end{example}

In what follows, we shall recall from {\cite[$\S$2.2]{Hoang}} how infinitesimal $\P$-bimodules can be represented as $\cS$-valued enriched functors. We will need to use the following notations and conventions.

\smallskip

$\bullet$ We write $\Fin_*$ to denote {the skeleton of} the \textbf{category of finite pointed sets}. The objects of $\Fin_*$ will be written as $\left \langle m \right \rangle := \{0,1,\cdots,m\}$  with $0$ as the base point, for $m\geqslant0$.

\smallskip

$\bullet$ A morphism $f: \l m \r \lrar \l n \r$ in $\Fin_*$ is called an \textbf{inert map} if $f^{-1}(i)$ is a singleton for every $i\in\{1,\cdots,n\}$, and an \textbf{active map} if $f^{-1}(0) = \{0\}$.

\smallskip

$\bullet$ We will write $\Fin_*^{\ine}$ (resp. $\Fin_*^{\act}$)  for the subcategory of $\Fin_*$ consisting of all the inert maps (resp. active maps).

\begin{conss}\label{infbimod} We construct three categories in $\cS$, denoted $\textbf{Ib}^{\mathcal{P}}$, $\textbf{R}^{\mathcal{P}}$ and  ${\textbf{I}\ell}^{\mathcal{P}}$, as follows. Firstly, the three categories have the same objects given by $C$-sequences $\{(c_1,\cdots,c_n;c) \, | \, c_i,c \in C , n\geqslant0  \}$. 

\smallskip
	
	$\bullet$ The mapping spaces of $\textbf{Ib}^{\mathcal{P}}$ are defined as follows. For each map $f : \l m \r \lrar \l n \r$ in $\Fin_*$, we define
	$$  \Map^f_{\textbf{Ib}^{\mathcal{P}}}((c_1,\cdots,c_n;c) , (d_1,\cdots,d_m;d)) := \mathcal{P}\left (c,\{d_j\}_{j\in f^{-1}(0)};d \right ) \otimes \bigotimes_{i=1,\cdots,n} \mathcal{P} \left (\{d_j\}_{j\in f^{-1}(i)};c_i \right ).$$
	Then we define
	$$ \Map_{\textbf{Ib}^{\mathcal{P}}}((c_1,\cdots,c_n;c) , (d_1,\cdots,d_m;d)) := \bigsqcup_{\left \langle m \right \rangle  \overset{f }{\rar} \left \langle n \right \rangle} \Map^f_{\textbf{Ib}^{\mathcal{P}}}((c_1,\cdots,c_n;c) , (d_1,\cdots,d_m;d) ) $$
	where the coproduct ranges over $\Hom_{\Fin_*}(\l m \r, \l n \r)$. The unit morphisms of $\textbf{Ib}^{\mathcal{P}}$ are defined via the unit operations of $\mathcal{P}$, and moreover, the categorical structure maps are induced by the composition in $\mathcal{P}$. 
	
	\smallskip
	
	$\bullet$ The mapping spaces of $\textbf{R}^{\mathcal{P}}$ are concentrated in objects of the form
	$$ \Map_{\textbf{R}^{\mathcal{P}}}(\, (c_1,\cdots,c_n;c) , (d_1,\cdots,d_m;c) \,) := \bigsqcup_{\left \langle m \right \rangle  \underset{\text{active}}{\overset{f}{\lrar}} \left \langle n \right \rangle} \left [  \bigotimes_{i=1,\cdots,n} \mathcal{P} \left (\{d_j\}_{j\in f^{-1}(i)};c_i \right )  \right ] $$
	where the coproduct ranges over those maps $f\in\Hom_{\Fin_*^{\act}}(\l m \r,\l n \r)$. Observe that there is a map 
	$$ \Map_{\textbf{R}^{\mathcal{P}}}( \, (c_1,\cdots,c_n;c) , (d_1,\cdots,d_m;c) \, ) \lrar \Map_{\textbf{Ib}^{\mathcal{P}}}(\, (c_1,\cdots,c_n;c) , (d_1,\cdots,d_m;c) \,) $$
	induced by the embedding $\Fin_*^{\act} \lrar \Fin_*$ and by inserting the unit operation $\id_c$ into the factor $\P(c;c)$ of the right hand side. The categorical structure of $\textbf{R}^{\mathcal{P}}$ is then defined via the operad structure of $\P$, so that $\textbf{R}^{\mathcal{P}}$ forms a subcategory of $\textbf{Ib}^{\mathcal{P}}$.
	
	\smallskip
	
	$\bullet$ The mapping spaces of ${\textbf{I}\ell}^{\mathcal{P}}$ are given by
	$$  \Map_{\textbf{I}\ell^{\mathcal{P}}}((c_1,\cdots,c_n;c) , (d_1,\cdots,d_m;d)) := \bigsqcup_{f} \mathcal{P}\left (c,\{d_j\}_{j\in f^{-1}(0)};d \right )$$
	where the coproduct ranges over the subset of $\Hom_{\Fin_*^{\ine}}(\l m \r,\l n \r)$ consisting of those maps $f$ with the property that the two colors $c_i$ and $d_{f^{-1}(i)}$ coincide for every $i\in\{1,\cdots,n\}$. As above, there is an embedding $\textbf{I}\ell^{\mathcal{P}} \lrar \textbf{Ib}^{\mathcal{P}}$ induced by the canonical embedding $\Fin_*^{\ine} \lrar \Fin_*$ and by inserting the identity operation $\id_{c_i}$ into the factor $\P(c_i;c_i)$ for $i=1,\cdots,n$.
\end{conss}

\begin{prop}\label{p:ibfunc} There are natural isomorphisms of categories
	
	\smallskip
	
	(i)  $ \Fun(\textbf{Ib}^{\mathcal{P}}, \cS) \cong \IbMod(\P)$,
	
	\smallskip
	
	(ii)  $ \Fun(\textbf{R}^{\mathcal{P}}, \cS) \cong \RMod(\P)$, and
	
	\smallskip
	
	(iii)  $ \Fun(\textbf{I}\ell^{\mathcal{P}}, \cS) \cong \ILMod(\P)$.
	\begin{proof} The first two isomorphisms are included in {\cite[$\S$2.2]{Hoang}}. The other can be proved in the same manner.
	\end{proof}
\end{prop}

We provide additional examples and remark regarding Constructions \ref{infbimod}.

\begin{example}\label{ex:ibcom} It is clear that $\textbf{Ib}^{\Com}$ is isomorphic to $\Fin_*^{\op}$, and while $\textbf{I}\ell^{\Com}$ (resp. $\textbf{R}^{\Com}$) is isomorphic to $(\Fin_*^{\ine})^{\op}$ (resp. $(\Fin_*^{\act})^{\op}$). 
\end{example}

\begin{rem}\label{ex:ibnucom} Let $\P$ be such that $\P(c) = \emptyset_\cS$ for every color $c$. Then by construction, we have that $$\Map^f_{\textbf{Ib}^{\mathcal{P}}}((c_1,\cdots,c_n;c) , (d_1,\cdots,d_m;d)) = \emptyset_\cS$$ whenever $f$ is not a surjection. The same observation holds for $\textbf{I}\ell^{\P}$ and  $\textbf{R}^{\P}$. Let us now denote by $\nucom$ the \textbf{nonunital commutative operad} which coincides with $\Com$ except that $\nucom(0)=\emptyset_\cS$. Accordingly, we can see that $\textbf{Ib}^{\nucom} \cong (\Fin_*^{\sur})^{\op}$ where $\Fin_*^{\sur}$ denotes the subcategory of $\Fin_*$ consisting of surjective maps. We have also that
	$$ \textbf{I}\ell^{\nucom} \cong (\Fin_*^{\ine, \sur})^{\op} \;\; \text{and} \;\; \textbf{R}^{\nucom} \cong (\Fin_*^{\act, \sur})^{\op}$$
	in which $\Fin_*^{\ine, \sur} := \Fin_*^{\ine} \cap \Fin_*^{\sur}$ and $\Fin_*^{\act, \sur} := \Fin_*^{\act} \cap \Fin_*^{\sur}$. 
\end{rem}

\begin{example}\label{ex:ibass} By construction, the category $\textbf{Ib}^{\Ass}$ has the same objects as $\Fin_*$. The data of a morphism in $\Map_{\textbf{Ib}^{\Ass}}(\l n \r,\l m \r)$ consist of a map $f : \l m \r \lrar \l n \r$ in $\Fin_*$, and in addition, a tuple of permutations $$(\sigma_0,\sigma_1,\cdots,\sigma_n) \in \Sigma_{k_0}\times\Sigma_{k_1}\times\cdots\times\Sigma_{k_n}$$
	with $k_i$ being the cardinality of $f^{-1}(i)$. {Thus $\textbf{Ib}^{\Ass} \cong \Gamma(\as)^{\op}$ for  $\Gamma(\as)$ of \cite{Rich}}. Moreover, the significant functor $\bar{\sB} : \Gamma(\as)^{\op} \lrar \cS$ considered {therein} is nothing but a model for $\Ass^{\Inf} \in \IbMod(\Ass)$.
\end{example}

\section{Hochschild and  Quillen cohomologies of operadic algebras}\label{s:op}

In this section, we assume that the base category $\cS$ is {(at least)} \textbf{sufficient} in the sense of {\cite[$\S$3]{Hoang}}, and we let $\P$ be a $C$-colored operad in $\cS${, for some fixed set of colors $C$}. 

\smallskip

{We will define \textbf{spectral Hochschild cohomology} of operads and their algebras, and describe the passage between various categories of operadic bimodules and algebras. Furthermore, we demonstrate that the cohomologies of $\P$-algebras are governed by those of the operad $\P$ itself.}

\subsection{Spectral Hochschild cohomology}\label{s:spechoch}

Let us start with the notion of spectral Hochschild cohomology of operads.

\begin{dfn} The \textbf{Hochschild complex of $\P$} {is defined to be} the cotangent complex $\rL_{\P^{\Inf}} \in \T_{\P^{\Inf}}\IbMod(\P)$ of $\P^{\Inf} \in \IbMod(\P)$ (see $\S$\ref{s:not}). Moreover, the \textbf{(spectral) Hochschild cohomology} of $\P$ with coefficients in a given object $M\in\T_{\P^{\Inf}}\IbMod(\P)$ is the space
	$$ \HHH^\star(\P ; M) := \Map^{\h}_{\T_{\P^{\Inf}}\IbMod(\P)}(\rL_{\P^{\Inf}}, M).$$
	For each $n\in\ZZ$, the \textit{$n$-th Hochschild cohomology group of $\P$ with coefficients in $M$} is 
	$$ \HHH^n(\P  ; M) := \pi_0 \Map^{\der}_{\T_{\P^{\Inf}}\IbMod(\P)} (\rL_{\P^{\Inf}} , M[n]).$$
\end{dfn}

\begin{rem} By definition, the Hochschild cohomology of $\P$ coincides with Quillen cohomology of $\P$ considered as an infinitesimal bimodule over itself. 
\end{rem}

Let $A\in \Alg_\P(\cS)$ be an algebra over $\P$. The corresponding cohomology of $A$ is defined as follows.
\begin{dfn} {The \textbf{Hochschild complex of $A$} {is defined to be} the cotangent complex $\rL_{A^{\me}} \in \T_{A^{\me}} \Mod_\P^A$ of $A^{\me} \in \Mod_\P^A$. Moreover,} the \textbf{(spectral) Hochschild cohomology} of $A$ with coefficients in a given object $N\in\T_{A^{\me}} \Mod_\P^A$ is the space
	$$ \HHH^\star(A ; N) := \Map^{\h}_{\T_{A^{\me}} \Mod_\P^A}(\rL_{A^{\me}}, N).$$
	For each $n\in\ZZ$, the \textit{$n$-th Hochschild cohomology group of $A$ with coefficients in $N$} is given by
	$$ \HHH^n(A ; N) := \pi_0 \Map^{\der}_{\T_{A^{\me}} \Mod_\P^A}(\rL_{A^{\me}}, N[n]).$$
\end{dfn}

\begin{notns} By abuse of notation, {we let $\P\otimes\mathbb{S}\in\T_{\P^{\Inf}}\IbMod(\P)$ and $A\otimes\mathbb{S}\in\T_{A^{\me}} \Mod_\P^A$ denote the suspension spectra respectively given by
\begin{gather*} (\P\otimes\mathbb{S})_{n,n}(c_1,\cdots,c_m;c) := \P(c_1,\cdots,c_m;c)\otimes \sS^{n}, \; \text{and} \\ (A\otimes\mathbb{S})_{n,n}(c) := A(c)\otimes \sS^{n}
\end{gather*}		
		where $\sS^{n}$ refers to the $n$-sphere in $\cS$ (see Notations \ref{no:sphere}).}
\end{notns}

\begin{prop}\label{p:hhop} Suppose that $\P$ is levelwise cofibrant. Then $\P\otimes\,\mathbb{S}$ is a suspension spectrum model for the Hochschild complex $\rL_{\P^{\Inf}}$. Likewise, when $A$ is levelwise cofibrant then $A\otimes\mathbb{S}$ is a suspension spectrum model for $\rL_{A^{\me}}$.	
	\begin{proof} When $\P$ is levelwise cofibrant then the coproduct $\P^{\Inf} \sqcup \P^{\Inf}\in\IbMod(\P)$ has already the right homotopy type. Thus we may exhibit the constant spectrum $\Sigma^\infty(\P^{\Inf} \sqcup \P^{\Inf})$ as a model for $\rL_{\P^{\Inf}}\in\T_{\P^{\Inf}}\IbMod(\P)$. Moreover, observe that the $n$-suspension $\Sigma^n(\P^{\Inf} \sqcup \P^{\Inf})\in\IbMod(\P)_{\P^{\Inf}//\P^{\Inf}}$ is given on each level by
		$$ \Sigma^n(\P^{\Inf} \sqcup \P^{\Inf})(c_1,\cdots,c_m;c) \simeq \P(c_1,\cdots,c_m;c)\otimes \sS^{n}.$$
		This follows from the fact that colimits of infinitesimal $\P$-bimodules can be taken levelwise. We just showed that $\P\otimes\mathbb{S}$ is a suspension spectrum model for $\rL_{\P^{\Inf}}$ (see also Remark \ref{r:sus}). The second claim is very similar.  
	\end{proof}
\end{prop}

\begin{example} Recall from Remark \ref{r:encode} that there exists an operad $\textbf{O}_C$ such that the category of $\textbf{O}_C$-algebras is isomorphic to the category of $C$-colored operads in $\cS$. Thus we may consider $\P$ as an object of $\Alg_{\textbf{O}_C}(\cS)$. It can be proved further that there is a categorical isomorphism 
	$$  \Mod_{\textbf{O}_C}^\P  \cong \IbMod(\P) $$
	between $\P$-modules over $\textbf{O}_C$ and infinitesimal $\P$-bimodules, which identifies $\P^{\Inf}\in\IbMod(\P)$ with $\P$ itself as an object of $\Mod_{\textbf{O}_C}^\P$. Accordingly, we obtain that Hochschild cohomology of the operad $\P$ agrees with Hochschild cohomology of $\P\in\Alg_{\textbf{O}_C}(\cS)$. (On other hand, we can see that Quillen cohomology of $\P\in\Alg_{\textbf{O}_C}(\cS)$ coincides with the reduced Quillen cohomology of the operad $\P$).
\end{example}

\begin{rem} {Hochschild cohomology of the operad $\Com$ is trivial for any base category $\cS$}. This is because of the fact that $\Com^{\Inf} \in \IbMod(\Com)$ is a free object generated by $\mathcal{E}_*$ (see Example \ref{r:comfree}). An analogous observation holds for every commutative algebra $A$, due to the fact that $A^{\me} \in \Mod_{\Com}^A$ coincides with the free object generated by $1_\cS$.
\end{rem}

\begin{rem}\label{r:same1} Suppose that the pair $(\cS, \P)$ is sufficient, so that we obtain Quillen equivalences
	\begin{equation}\label{eq:seqQuillen1}
		\adjunction*{\simeq}{\T_{\P^{\Inf}}\IbMod(\P)}{\T_{\P^{\sB}}\BMod(\P)}{} \adjunction*{\simeq}{}{\T_\P\Op_C(\mathcal{S})}{} \adjunction*{\simeq}{}{\T_\P\Op(\mathcal{S})}{}
	\end{equation}
	(see $\S$\ref{s:rv}). In particular, this allows us to calculate both Quillen and Hochschild cohomologies of $\P$ within the same category $\T_{\P^{\Inf}}\IbMod(\P)$, and furthermore, the same category $\IbMod(\P)$ when $(\cS, \P)$ is stably sufficient. 
\end{rem}

\begin{rem}\label{r:same2} Similarly, when the triple $(\cS,\P,A)$ is sufficient, then due to the Quillen equivalence
	$$ \adjunction*{\simeq}{\T_{A^{\me}} \Mod_\P^A}{\T_A\Alg_\P(\cS)}{},$$
	we may also consider both Quillen and Hochschild cohomologies of $A$ in the same framework of $\T_{A^{\me}} \Mod_\P^A$. Moreover, when $(\cS,\P,A)$ is stably sufficient then the two cohomologies can be brought into the same category $\Mod_\P^A$.
\end{rem}

\begin{example}\label{ex:hhpstab} {When} the pair $(\cS,\P)$ is stably sufficient then under the right Quillen equivalence $$ \ker\circ \, \Omega^\infty : \T_{\P^{\Inf}}\IbMod(\P) \lrarsimeq \IbMod(\P),$$ the object $\rL_{\P^{\Inf}}$ is simply identified with $\P^{\Inf}\in\IbMod(\P)$ (cf. {\cite[Corollary 2.2.4]{Yonatan}}). Therefore, after sending coefficients into $\IbMod(\P)$, the Hochschild cohomology of $\P$ with coefficients in {an object} $M\in\IbMod(\P)$ is given by
	$$ \HHH^\star(\P ; M) \simeq \Map^{\h}_{\IbMod(\P)}(\P^{\Inf}, M).$$
	In particular, {when} $\P$ is a dg operad concentrated in arity $1$ (i.e. $\P$ is identified with a dg category), then we recover the usual \textbf{Hochschild cohomology of a dg category}, that has been widely considered in literature {(cf., e.g., \cite{Toen})}.
\end{example}

\begin{example}\label{ex:hhastab} In the same fashion, when the triple $(\cS,\P, A)$ is stably sufficient,  the right Quillen equivalence $ \ker\circ \, \Omega^\infty : \T_{A^{\me}}\Mod_\P^A \lrarsimeq \Mod_\P^A$ identifies $\rL_{A^{\me}} \in \T_{A^{\me}}\Mod_\P^A$ to $A^{\me}\in\Mod_\P^A$. Thus in this situation, the Hochschild cohomology of $A$ with coefficients in a given object $N\in\Mod_\P^A$ is computed by
	$$ \HHH^\star(A ; N) \simeq \Map^{\h}_{\Mod_\P^A}(A^{\me}, N).$$
\end{example}

\begin{example} For a dg associative algebra $A$, {the category $\Mod_{\Ass}^A$ coincides with that} of $A$-\textbf{bimodules}, and hence, in this situation we recover the classical Hochschild cohomology of dg associative algebras. 
\end{example}

The following example will be clarified in $\S$\ref{s:cotan}.

\begin{example} {For a fibrant simplicial operad $\P$}, we have an equivalence of $\infty$-categories
	$$\T_{\P^{\Inf}}\IbMod(\P)_\infty \simeq \Fun(\Tw(\P), \Sp)$$
	in which $\Tw(\P)$ is the \textit{twisted arrow $\infty$-category of $\P$} and $\Sp$ denotes the $\infty$-category of spectra. Moreover, it can be shown that under this equivalence, $\rL_{\P^{\Inf}} \in \T_{\P^{\Inf}}\IbMod(\P)$ is identified with the constant functor on the sphere spectrum $\mathbb{S}\in\Sp$. Thus after sending coefficients into $\Fun(\Tw(\P), \Sp)$, the Hochschild cohomology of $\P$ with coefficients in a given functor $\F:\Tw(\P) \lrar \Sp$ is simply classified by  
	$\lim\F$. In particular, when $\P$ is concentrated in arity $1$ (i.e. $\P$ is identified with a simplicial category or an $\infty$-category), then we obtain the usual cohomology of $\infty$-categories which generalizes the \textbf{Baues-Wirsching cohomology} of small categories (cf. \cite{BauWi}).
\end{example}

\subsection{Endomorphism constructions}\label{s:endo}

It is known that the composite product $$ - \circ - : \Coll_C(\cS) \times \Coll_C(\cS) \lrar \Coll_C(\cS)$$ is \textit{left linear}, i.e. for every $A \in \Coll_C(\cS)$ the functor $ (-) \circ A : \Coll_C(\cS) \lrar \Coll_C(\cS)$ preserves all colimits. Consider the case where $A \in \cS^{\times C}$ regarded as a $C$-collection concentrated in level $0$. Taking composite product with $A$ determines a functor $$ (-) \circ A : \Coll_C(\cS) \lrar  \cS^{\times C}.$$
Of course, this functor admits a right adjoint, as described in the following construction.

\begin{cons} Let $A=\{A(c)\}_{c\in C}$ and $B=\{B(c)\}_{c\in C}$ be two objects of $\cS^{\times C}$. The \textbf{endomorphism object} associated to the  pair $(A,B)$ is a $C$-collection denoted by $\End_{A,B}$, and defined by letting
	$$ \End_{A,B}(c_1,\cdots,c_n;c) := \Map_{\cS}(A(c_1)\otimes\cdots\otimes A(c_n),B(c)).$$
	The right $\Sigma_n$-action is naturally given by the permutation on the factors $A(c_i)$'s. In particular, the \textbf{endomorphism operad associated to $A$} is given by $\End_A:=\End_{A,A}$. (See also \cite{Fresse1, Rezk}).
\end{cons}

We have indeed the following result, which can be readily verified using definition.

\begin{prop} The functors $ (-) \circ A$ and $\End_{A,-}$ determine an adjunction
	\begin{equation}\label{eq:endbasic}
		(-) \circ A : \adjunction*{}{\Coll_C(\cS)}{\cS^{\times C}}{} : \End_{A,-}.
	\end{equation}
\end{prop}

One may obtain relative versions of \eqref{eq:endbasic} when endowing the involving objects with additional structures, as discussed in what follows.

\smallskip

Let us assume further that $A$ comes equipped with a $\P$-algebra structure classified by a map $\ell_A : \P \lrar \End_A$ of operads. The \textbf{relative composite product} determines a functor $$(-)\circ_\P A : \RMod(\P)\lrar \cS^{\times C}.$$
On other hand, for any $B \in \cS^{\times C}$, one observes that $\End_{A,B}$ carries a canonical right module structure over $\End_A$.  Thus, $\End_{A,B}$ carries a right $\P$-module structure given by the restriction along $\ell_A$. Furthermore, when $B$ comes equipped with the structure of a $\P$-algebra then $\End_{A,B}$ carries a canonical $\P$-bimodule structure. The relative versions of \eqref{eq:endbasic} are as follows.

\begin{prop}\cite{Fresse1, Rezk} \label{p:end} {Let $A$ be a $\P$-algebra. Then} the constructions $(-)\circ_\P A$ and  $\End_{A,-}$ determine the adjunctions:	
	\begin{equation}\label{eq:end1}
		(-)\circ_\P A : \adjunction*{}{\RMod(\P)}{\cS^{\times C}}{} : \End_{A,-}, \;\; \text{and}
	\end{equation}
	\begin{equation}\label{eq:end2}
		(-)\circ_\P A : \adjunction*{}{\BMod(\P)}{\Alg_\P(\cS)}{} : \End_{A,-}.
	\end{equation}
\end{prop}

Another important adjunction, which interposes the two above, is given as follows.

\begin{prop}\label{p:end3} {Let $A$ be a $\P$-algebra. Then} the constructions $(-)\circ_\P A$ and  $\End_{A,-}$ determine an adjunction:	
	\begin{equation}\label{eq:end3}
		(-)\circ_\P A : \adjunction*{}{\IbMod(\P)}{\Mod_\P^A}{} : \End_{A,-}.
	\end{equation}
	\begin{proof} {For an $A$-module $N$, the canonical infinitesimal $\P$-bimodule structure on $\End_{A,N}$ is defined as follows. Since $\End_{A,N}$ already inherits a right $\P$-module structure, it suffices to establish a compatible infinitesimal left $\P$-module structure, whose data consist of $\Sigma_*$-equivariant maps of the form
			$$ \P(c_1,\cdots,c_n;c) \otimes \End_{A,N}(d_1,\cdots,d_m;c_i) \lrar \End_{A,N}(c_1,\cdots,c_{i-1},d_1,\cdots,d_m,c_{i+1},\cdots,c_n;c)$$ 
			(see Definition \ref{d:infmodules}). This is equivalent to giving a map of the form
			$$ \P(c_1,\cdots,c_n;c)  \otimes \textbf{A}^{\otimes\hat{i}} \otimes [A(d_1) \otimes \cdots \otimes A(d_m)] \otimes \Map_{\cS}(A(d_1) \otimes \cdots \otimes A(d_m),N(c_i)) \lrar N(c) $$
			where $\textbf{A}^{\otimes\hat{i}} := \bigotimes_{k\in \{1,\cdots,n\} \setminus \{i\}} A(c_k)$. This map is defined by first inserting the canonical map $$ [A(d_1) \otimes \cdots \otimes A(d_m)] \otimes \Map_{\cS}(A(d_1) \otimes \cdots \otimes A(d_m),N(c_i)) \lrar N(c_i),$$
			and then applying the structure map $\P(c_1,\cdots,c_n;c)\otimes \textbf{A}^{\otimes\hat{i}} \otimes N(c_i) \lrar N(c)$ of $N$. The axioms of an infinitesimal $\P$-bimodule structure follow from those of the $A$-module structure on $N$.}  
		
		\smallskip
		
		{Next, the object  $M\circ_\P A$ for some $M \in \IbMod(\P)$ carries a canonical $A$-module structure, described as follows. By definition, we need to establish $\Sigma_*$-equivariant maps of the form
			$$  \P(c_1,\cdots,c_n;c) \otimes \textbf{A}^{\otimes\hat{i}} \otimes (M\circ_\P A)(c_i) \lrar (M\circ_\P A)(c).$$
			This can be constructed from maps of the form
			$$ \P(c_1,\cdots,c_n;c) \otimes \textbf{A}^{\otimes\hat{i}} \otimes M(d_1,\cdots,d_m;c_i) \otimes [A(d_1) \otimes \cdots \otimes A(d_m)] $$
			$$ \lrar M(c_1,\cdots,c_{i-1},d_1,\cdots,d_m,c_{i+1},\cdots,c_n;c) \otimes \textbf{A}^{\otimes\hat{i}} \otimes [A(d_1) \otimes \cdots \otimes A(d_m)] $$
			that satisfy a compatibility condition governed by the coequalizer $M\circ\P\circ A \mathrel{\substack{\longrightarrow \\ \longrightarrow}} M \circ A$. This is simply induced by the structure map of $M$:
			$$  \P(c_1,\cdots,c_n;c) \otimes  M(d_1,\cdots,d_m;c_i)  \lrar M(c_1,\cdots,c_{i-1},d_1,\cdots,d_m,c_{i+1},\cdots,c_n;c).$$}   
		
		{Finally, to verify that the obtained functors form an adjoint pair, we use the adjunction \eqref{eq:end1}, reducing the task to showing that the unit $M \lrar \End_{A,M\circ_\P A}$ (resp. counit $\End_{A,N}\circ_\P A \lrar N$) of \eqref{eq:end1} is a map of infinitesimal $\P$-bimodules (resp. $A$-modules) for $M \in \IbMod(\P)$ (resp. $N \in \Mod_\P^A$). This follows straightforwardly from the definitions.}
	\end{proof}
\end{prop}

A noteworthy consequence is as follows.

\begin{cor}\label{r:freeA} For an object $X\in\cS^{\times C} \subseteq \Coll_C(\cS)$, there is an isomorphism of $A$-modules over $\P$:
	$$ \Free_A(X) \cong (\mathcal{P}\circ_{(1)}X) \circ_\P A.$$
	In particular, when $C = \{*\}$ we obtain that $\Free_A(1_\cS) \cong  \mathcal{P}[1]\circ_\P A$ (see $\S$\ref{s:infbi} for notations).
	\begin{proof} Consider the following square of left adjoints
		$$ \xymatrix{
			\Coll_C(\cS) \ar[r]^{\,\; (-) \circ A}\ar[d]_{\Free_\P^{\Inf}} & \cS^{\times C} \ar[d]^{\Free_A} \\
			\IbMod(\P) \ar[r]_{\,\; (-) \circ_\P A} & \Mod_\P^A \\
		}$$
		which is commutative due to that property of the associated square of right adjoints. Combined with \eqref{eq:freeinf}, this yields for each $X\in\cS^{\times C} \subseteq \Coll_C(\cS)$ a chain of isomorphisms in $\Mod_\P^A$:
		$$ \Free_A(X) \cong \Free_A(X \circ A) \cong \Free_\P^{\Inf}(X)\circ_\P A \cong (\mathcal{P}\circ_{(1)}X) \circ_\P A.$$
		For the case $C=\{*\}$, we obtain that $$\Free_A(1_\cS) \cong (\mathcal{P}\circ_{(1)}\mathcal{E}_*)\circ_\P A \cong \mathcal{P}[1]\circ_\P A.$$  
	\end{proof}
\end{cor}

\begin{rem} The latter identification can also be found in {\cite[$\S$4.3]{Fresse1}}. Due to this, one may obtain a description of $\Free_A(1_\cS)$, which models the universal enveloping algebra of $A$ (see Proposition \ref{p:uenv}).
\end{rem}

\begin{example} An interesting example concerns the two operads $\Lie$ and $\Ass$ enriched over $\cS = \C(\textbf{k})$ (see Examples \ref{ex:base}). The canonical embedding $\Lie \lrar \Ass$ endows $\Ass$ with an infinitesimal $\Lie$-bimodule structure. Moreover, it can be shown that the map $\Free_{\Lie}^{\Inf}(\mathcal{E}_*) = \Lie[1]  \overset{\cong}{\lrar} \Ass$ classified by the unique nullary operation of $\Ass$ is an isomorphism in $\IbMod(\Lie)$. Combined with Corollary \ref{r:freeA}, this yields for each $\Lie$-algebra $\g$ a chain of natural isomorphisms in $\Mod^\g_{\Lie}$:
	$$ \Ass\circ_{\Lie}\g \cong \Free_{\Lie}^{\Inf}(\mathcal{E}_*)\circ_{\Lie}\g \cong \Free_\g(\textbf{k}) \cong \sU_{\Lie}(\g).$$
	This allows us to recover a fundamental result asserting that the universal enveloping algebra of $\g$ can be naturally modeled by $\Ass\circ_{\Lie}\g$, i.e. the induction of $\g$ along the map $\Lie \lrar \Ass$.
\end{example}

To be able to do homotopy theory, we will need the following assertion.

\begin{prop} The adjunctions \eqref{eq:end1}, \eqref{eq:end2} and \eqref{eq:end3} are all Quillen adjunctions, provided that $A\in\Alg_\P(\cS)$ is levelwise cofibrant. 
	\begin{proof} Clearly with the assumption that $A$ is levelwise cofibrant, the construction $\End_{A,-}$ preserves fibrations and trivial fibrations, due to basic properties of a symmetric monoidal model category.
	\end{proof}	
\end{prop}

\begin{rem}\label{r:keydiagram} By the identification $\P\circ_\P A \cong A$, we obtain further a commutative diagram of adjunctions: 
	\begin{equation}\label{eq:prekey}
		\begin{tikzcd}[row sep=3.5em, column sep =3.5em]
			\RMod(\P)_{\P^{\sr}/}	 \arrow[r, shift left=1.5] \arrow[d, shift right=1.5] & \IbMod(\P)_{\P^{\Inf}/} \arrow[l, shift left=1.5, "\perp"'] \arrow[d, shift left=1.5] \arrow[r, shift left=1.5]  & \BMod(\P)_{\P^{\sB}/}  \arrow[d, shift left=1.5] \arrow[l, shift left=1.5, "\perp"']   \\
			(\cS^{\times C})_{A/} \arrow[r, shift right=1.5] \arrow[u, shift right=1.5, "\dashv"] & (\Mod_\P^A)_{A^{\me}/}  \arrow[l, shift right=1.5, "\downvdash"] \arrow[u, shift left=1.5, "\vdash"'] \arrow[r, shift right=1.5] & (\Alg_\P(\cS))_{A/}  \arrow[u, shift left=1.5, "\vdash"'] \arrow[l, shift right=1.5, "\downvdash"]
		\end{tikzcd}
	\end{equation}
	in which the vertical adjunctions are all induced by the adjunction $(-)\circ_\P A \dashv \End_{A,-}$, and the horizontal pairs are given by the induction-restriction adjunctions.
\end{rem}

\subsection{Formulations and comparison theorems}\label{s:op2}

We are interested in the right square of \eqref{eq:prekey}, from which we obtain a commutative square of Quillen adjunctions between tangent categories:
\begin{equation}\label{eq:key}
	\begin{tikzcd}[row sep=3.5em, column sep =3.5em]
		\mathcal{T}_{\P^{\Inf}}\IbMod(\P) \arrow[r, shift left=1.5] \arrow[d, shift right=1.5] & \mathcal{T}_{\P^{\sB}}\BMod(\P) \arrow[l, shift left=1.5, "\perp"'] \arrow[d, shift left=1.5]  & \\
		\T_{A^{\me}} \Mod_\P^A \arrow[r, shift right=1.5] \arrow[u, shift right=1.5, "\dashv"] & \T_A\Alg_\P(\cS)  \arrow[l, shift right=1.5, "\downvdash"] \arrow[u, shift left=1.5, "\vdash"'] &
	\end{tikzcd}
\end{equation}

\begin{notn}\label{no:circend} We will use the same notation $(-)\circ^{\st}_\P A \dashv \End^{\st}_{A,-}$ to signify the two vertical adjunctions of \eqref{eq:key}:
	\begin{gather*}
		(-)\circ^{\st}_\P A : \adjunction*{}{\mathcal{T}_{\P^{\Inf}}\IbMod(\P)}{\T_{A^{\me}} \Mod_\P^A}{} : \End^{\st}_{A,-} \;\;\; \text{and} \\ \;\;  (-)\circ^{\st}_\P A : \adjunction*{}{\mathcal{T}_{\P^{\sB}}\BMod(\P)}{\T_A\Alg_\P(\cS)}{} : \End^{\st}_{A,-}.
	\end{gather*}
\end{notn}

\begin{rem}\label{r:endcons} By construction, for each $M\in\mathcal{T}_{\P^{\sB}}\BMod(\P)$ (or $\mathcal{T}_{\P^{\Inf}}\IbMod(\P)$), the spectrum object $M\circ^{\st}_\P A$ is simply taken degreewise, i.e.  $$(M\circ^{\st}_\P A)_{m,n}=M_{m,n}\circ_\P A.$$ 
	The case of $\End^{\st}_{A,-}$ is rather different. Namely, for each $N\in\T_A\Alg_\P(\cS)$ (or $\T_{A^{\me}}\Mod_\P^A$), the spectrum object $\End^{\st}_{A,N}$ is given at each bidegree $(m,n)$ as the pullback
	$$ \xymatrix{
		(\End^{\st}_{A,N})_{m,n} \ar[r]\ar[d] & \End_{A,N_{m,n}} \ar[d] \\
		\P \ar[r]^{\ell_A} & \End_A \\
	}$$
	where the map $\End_{A,N_{m,n}} \lrar \End_A$ ($\x{\defi}{=} \End_{A,A}$) is induced by the structure map $N_{m,n} \lrar A$.
\end{rem}

The following plays a pivotal role in the study of various cohomologies of operadic algebras.
\begin{prop}\label{p:bialg}  Suppose that $A\in\Alg_\P(\cS)$ is levelwise cofibrant. There are weak equivalences $$\rL_{\P^{\sB}}\circ^{\st}_\P A \simeq \rL_A \;\; \text{and} \;\; \rL_{\P^{\Inf}}\circ^{\st}_\P A \simeq \rL_{A^{\me}}$$
	in $\T_A\Alg_\P(\cS)$ and $\T_{A^{\me}} \Mod_\P^A$ respectively.
	\begin{proof} We will prove the first weak equivalence. The other can then be verified in the same fashion. By definition, $\rL_{\P^{\sB}}$ is represented by a constant spectrum of the form $\Sigma^\infty\left(\P^{\sB}\sqcup(\P^{\sB})^{\cof}\right) \in \mathcal{T}_{\P^{\sB}}\BMod(\P)$ in which $(\P^{\sB})^{\cof}$ refers to any cofibrant model for $\P^{\sB} \in \BMod(\P)$, and the coproduct $\P^{\sB}\sqcup(\P^{\sB})^{\cof}$ is considered as an object over and under $\P^{\sB}$. We have a chain of weak equivalences
		$$ \rL_{\P^{\sB}}\circ^{\st}_\P A \simeq \Sigma^\infty\left(\P^{\sB}\sqcup(\P^{\sB})^{\cof}\right)\circ^{\st}_\P A \; \x{\defi}{=} \; \Sigma^\infty\left((\P^{\sB}\sqcup(\P^{\sB})^{\cof})\circ_\P A\right) \simeq \Sigma^\infty(A\sqcup A^{\cof})$$
		where $A^{\cof} := (\P^{\sB})^{\cof}\circ_\P A$ represents a cofibrant model for $A\in\Alg_\P(\cS)$. As above, $\Sigma^\infty(A\sqcup A^{\cof})$ is a model for $\rL_A$. So we obtain a weak equivalence $\rL_{\P^{\sB}}\circ^{\st}_\P A \simeq \rL_A$ as expected. 
	\end{proof}
\end{prop}

According to $\S$\ref{s:rv}, $\rL_{\P^{\sB}}\in\mathcal{T}_{\P^{\sB}}\BMod(\P)$ admits a suspension spectrum replacement given by $\P\circ\sS_C^{\bullet}$ with $(\P\circ\sS_C^{\bullet})_{n,n}=\P\circ\sS_C^{n}$,  provided that $\P$ is $\Sigma$-cofibrant. Due to this, we can describe the cotangent complex $\rL_A$ as follows.

\begin{prop}\label{p:cotalg}  Suppose that the triple $(\cS,\P, A)$ is sufficient. Then for every $n\in\mathbb{N}$ we have a weak equivalence in $\Alg_\P(\cS)_{A//A}$:
	$$ (\P\circ\sS_C^{n})\circ_\P A \simeq \Sigma^n(A \sqcup A) $$
	in which the right hand side is the $n$-suspension of $A\sqcup A \in \Alg_\P(\cS)_{A//A}$. Consequently, the cotangent complex $\rL_A \in \T_A\Alg_\P(\cS)$ admits a suspension spectrum replacement:
	$$ (\P\circ\sS_C^{\bullet})\circ^{\st}_\P A \simeq \rL_A.$$
	\begin{proof} As discussed in $\S$\ref{s:rv}, we have a weak equivalence $\P\circ\sS_C^{n} \simeq \Sigma^n(\P^{\sB} \sqcup \P^{\sB})$ in $\BMod(\P)_{\P^{\sB}//\P^{\sB}}$. By the assumption that $\P$ is $\Sigma$-cofibrant and that $A$ is cofibrant, the relative composite product $(\P\circ\sS_C^{n})\circ_\P A$ has the right homotopy type (cf. \cite{Fresse1}, $\S$15.2). Thus we obtain a chain of natural weak equivalences
		$$ (\P\circ\sS_C^{n})\circ_\P A \simeq \left(\Sigma^n(\P^{\sB} \sqcup \P^{\sB})\right) \circ_\P A \simeq \Sigma^n(A\sqcup A)$$
		where the latter weak equivalence is due to the left Quillen functor $$(-)\circ_\P A : \BMod(\P)_{\P^{\sB}//\P^{\sB}} \lrar \Alg_\P(\cS)_{A//A}.$$ The second statement follows from the above, because the functor $(-)\circ^{\st}_\P A$ is computed degreewise. (Of course, it can also be immediately derived from Proposition \ref{p:bialg}).
	\end{proof}
\end{prop}

In summary, we obtain the following conclusion.

\begin{thm}\label{t:maincotalg} Suppose that the triple $(\cS,\P, A)$ is sufficient. The Quillen cohomology of $A$ with coefficients in a given object $N\in\T_{A^{\me}} \Mod_\P^A$ is computed by the formula
	$$ \HHQ^\star(A ; N) \simeq \Map^{\h}_{\T_{A^{\me}} \Mod_\P^A}\left((\P\circ\sS_C^{\bullet})\circ^{\st}_\P A, N \right).$$
	When $(\cS,\P, A)$ is in addition stably sufficient, the Quillen cohomology of $A$ with coefficients in some $N\in\Mod_\P^A$ is computed by the formula
	$$ \HHQ^\star(A ; N) \simeq \Map^{\h}_{\Mod_\P^A}(\ovl{\rL}_\P\circ_\P A, N)$$
	{where the object $\ovl{\rL}_\P \in\IbMod(\P)$ is described in Theorem \ref{t:maincot}.}
	\begin{proof} First, under the right Quillen equivalence $\T_A\Alg_\P(\cS) \lrarsimeq \T_{A^{\me}} \Mod_\P^A$
		(cf. Theorem \ref{t:tanalg}), the object $(\P\circ\sS_C^{\bullet})\circ^{\st}_\P A$ is identified with itself $(\P\circ\sS_C^{\bullet})\circ^{\st}_\P A$ considered as a prespectrum of the underlying $A$-modules. (This follows from {\cite[Corollary 2.4.8]{YonatanBundle}}). Thus, the latter is a model for the derived image of $\rL_A$ in $\T_{A^{\me}} \Mod_\P^A$. This proves the first statement.
		
		\smallskip
		
		{Next, when $(\cS,\P, A)$ is stably sufficient, we obtain an extended commutative diagram of Quillen adjunctions of the form} 
		\begin{equation}\label{eq:key1}
			\begin{tikzcd}[row sep=3.5em, column sep =3.5em]
				\IbMod(\P)	 \arrow[r, shift left=1.5, "\simeq"] \arrow[d, shift right=1.5, "(-)\circ_\P A"'] & \mathcal{T}_{\P^{\Inf}}\IbMod(\P) \arrow[l, shift left=1.5, "\perp"'] \arrow[d, shift left=1.5, "(-)\circ^{\st}_\P A"] \arrow[r, shift left=1.5, "\simeq"]  & \mathcal{T}_{\P^{\sB}}\BMod(\P)  \arrow[d, shift left=1.5, "(-)\circ^{\st}_\P A"] \arrow[l, shift left=1.5,  "\perp"']   \\
				\Mod_\P^A \arrow[r, shift right=1.5] \arrow[u, shift right=1.5, "\dashv", "\End_{A,-}"'] & \T_{A^{\me}} \Mod_\P^A  \arrow[l, shift right=1.5, "\downvdash", "\simeq" '] \arrow[u, shift left=1.5, "\vdash"', "\End_{A,-}^{\st}"] \arrow[r, shift right=1.5] & \T_A\Alg_\P(\cS)  \arrow[u, shift left=1.5, "\vdash"', "\End_{A,-}^{\st}"] \arrow[l, shift right=1.5, "\downvdash", "\simeq" ']
			\end{tikzcd}
		\end{equation} 
		{in which the horizontal pairs are all Quillen equivalences. Recall from $\S$\ref{s:rv} that under the right Quillen equivalence $\mathcal{T}_{\P^{\Inf}}\IbMod(\P)  \lrarsimeq  \IbMod(\P)$, the object $\P\circ\sS_C^{\bullet}\in\mathcal{T}_{\P^{\Inf}}\IbMod(\P)$ is identified with $\ovl{\rL}_\P \in \IbMod(\P)$. Combined with the first paragraph and with the fact that $\ovl{\rL}_\P$ is $\Sigma$-cofibrant, it yields a weak equivalence in $\Mod_\P^A$ of the form $\ovl{\rL}_\P\circ_\P A \simeq \rL_A$. (Here we use the same notation to denote the derived image of $\rL_A \in \T_A\Alg_\P(\cS)$ in $\Mod_\P^A$.) This verifies the second statement.}  
	\end{proof}
\end{thm}

\begin{example} {For $\cS=\C(\textbf{k})$ with $\textbf{k}$ being a commutative ring (see Examples \ref{ex:base})}, we can prove that $\ovl{\rL}_\P\in\IbMod(\P)$ is weakly equivalent to the \textit{infinitesimal composite product} $\P\circ_{(1)}\I_C$ (see $\S$\ref{s:operad} for notation). Furthermore, it can be shown that the object $$\ovl{\rL}_\P\circ_\P A \simeq (\P\circ_{(1)}\I_C)\circ_\P A$$ is a model for the \textbf{module of K{\"a}hler differentials} $\Omega_A\in\Mod_\P^A$. {This allows us to recover a classical result} asserting that, in a differential graded setting, the cotangent complex of an operadic algebra is equivalent to its module of K{\"a}hler differentials (cf., e.g. \cite{YonatanCotangent, Loday, Fresse1}). Detailed discussion regarding this topic will be included in a subsequent paper.
\end{example}

In what follows, we shall provide some comparison theorems asserting that Hochschild and Quillen cohomologies of $A$ can be encoded by the corresponding cohomologies of $\P$ itself. 

\begin{thm} \label{co:bialg} {Suppose that the triple $(\cS,\P,A)$ is sufficient and} suppose given a fibrant object $N\in \T_{A^{\me}} \Mod_\P^A$. 
	
	\smallskip
	(i) There is a natural weak equivalence
	$$ \HHQ^\star(A ; N)  \simeq  \HHQ^\star(\P^{\sB} ; \End^{\st}_{A,N})$$
	between Quillen cohomology of $A$ with coefficients in $N$ and Quillen cohomology of $\P^{\sB}\in\BMod(\P)$ with coefficients in $\End^{\st}_{A,N} \in \mathcal{T}_{\P^{\Inf}}\IbMod(\P)$.
	
	\smallskip
	(ii) There is a natural weak equivalence
	$$ \HHH^\star(A ; N)  \simeq \HHH^\star(\P ; \End^{\st}_{A,N}) $$
	between the Hochschild cohomology of $A$ with coefficients in $N$ and Hochschild cohomology of $\P$ with coefficients in $\End^{\st}_{A,N}$. 
	
	\begin{proof} {First, by assumption the horizontal adjunctions of \eqref{eq:key} are all Quillen equivalences.}  After sending coefficients into $\mathcal{T}_{\P^{\Inf}}\IbMod(\P)$, we {get} that
		$$ \HHQ^\star(\P^{\sB} ; \End^{\st}_{A,N}) \simeq \Map^{\h}_{\mathcal{T}_{\P^{\Inf}}\IbMod(\P)}(\P\circ\sS_C^{\bullet},\End^{\st}_{A,N}),$$
		and while by Theorem \ref{t:maincotalg} we have $\HHQ^\star(A;N) \simeq \Map^{\h}_{\T_{A^{\me}} \Mod_\P^A}\left((\P\circ\sS_C^{\bullet})\circ^{\st}_\P A,N \right)$. Combined with the adjunction $(-)\circ^{\st}_\P A\dashv\End^{\st}_{A,-}$, these indeed verify (i). The statement (ii) is verified by combining the adjunction $(-)\circ^{\st}_\P A\dashv\End^{\st}_{A,-}$ with Proposition \ref{p:bialg}.
	\end{proof}
\end{thm}

We {now} arrive at the main result as follows.
\begin{thm}\label{t:opmain} {Suppose that the triple $(\cS,\P,A)$ is abundant.} For a fibrant object $N\in \T_{A^{\me}} \Mod_\P^A$, there is a natural weak equivalence
	$$\HHQ^\star(A;N)  \simeq \Omega\HHQ^\star(\P ; \End^{\st}_{A,N})$$
	where the loop space on the right is taken at the zero element (see Remark \ref{r:trivclass}). In particular, for every integer $n$ there is an isomorphism
	$$ \HHQ^n(A;N)  \cong \HHQ^{n+1}(\P;\End^{\st}_{A,N}).$$
	\begin{proof} We have in fact a chain of weak equivalences
		$$ \HHQ^\star(A ; N)  \simeq \HHQ^\star(\P^{\sB} ; \End^{\st}_{A,N}) \simeq \Omega\HHQ^\star(\P ; \End^{\st}_{A,N})$$
		where the first weak equivalence is due to Theorem \ref{co:bialg}(i), and the second weak equivalence follows from Proposition \ref{p:key}.
	\end{proof}
\end{thm}

\begin{rem} The above theorem admits a stable analogue as follows. Assume that the triple $(\cS,\P,A)$ is stably abundant, so that the horizontal adjunctions of \eqref{eq:key1} are all Quillen equivalences. Similarly as above, we may obtain a chain of weak equivalences 
	$$\HHQ^\star(A ; N)  \simeq  \HHQ^\star(\P^{\sB} ; \End_{A,N}) \simeq \Omega\HHQ^\star(\P ; \End_{A,N})$$
	in which $N$ is now a fibrant object in $\Mod_\P^A$, and $\End_{A,N}$ is simply an infinitesimal $\P$-bimodule. 
\end{rem}

\begin{rem} When the triple $(\cS,\P,A)$ is stably sufficient, an analogue for Hochschild cohomology is now evident, as described via the following chain of weak equivalences:
	$$\HHH^\star(A;N) \simeq \Map^{\h}_{\Mod_\P^A}(A^{\me},N) \simeq \Map^{\h}_{\IbMod(\P)}(\P^{\Inf},\End_{A,N}) \simeq \HHH^\star(\P ; \End_{A,N})$$
	(see Examples \ref{ex:hhpstab} and \ref{ex:hhastab}).
\end{rem}

\section{Simplicial operads and their associated algebras}\label{s:Qprin}

{In this section, the category $\Set_\Delta$ has the Kan-Quillen model structure,} and  $\Op(\Set_\Delta)$ is endowed with the Dwyer-Kan model structure (see $\S$\ref{s:optransfermod}). We let $\P$ be a fixed $C$-colored simplicial operad for some set $C$. 

\smallskip

{We first recall from \cite{Hoang} the construction of the twisted arrow $\infty$-category $\Tw(\P)$}, and recall how the tangent categories at $\P$ can be represented via that construction. Due to this, we may then formulate Hochschild and Quillen cohomologies of $\P$, as well as cohomologies of $\P$-algebras. In the last subsection, we shall formulate the Quillen principle for $\E_n$-operads, and from which we obtain the Quillen principle for algebras over $\E_n$.

\subsection{Operadic twisted arrow $\infty$-categories and tangent categories}\label{s:twis}

We will assume that $\P$ is fibrant (i.e. every space of operations of $\P$ is a Kan complex). Recall from \cite{Hoang} that the twisted arrow $\infty$-category $\Tw(\P)$ is given by the (covariant) unstraightening of the simplicial copresheaf $\P^{\Inf} : \textbf{Ib}^{\mathcal{P}} \lrar \Set_\Delta$ which encodes the datum of $\P$ as an infinitesimal bimodule over itself (see $\S$\ref{s:infbi}). In particular, $\Tw(\P)$ is endowed with a left fibration $\Tw(\P) \lrar \sN(\textbf{Ib}^{\mathcal{P}})$
where $\sN(-)$ denotes the \textit{simplicial nerve functor}.

\begin{rem} Unwinding the definitions, objects of $\Tw(\P)$ are precisely the \textbf{operations} of $\P$ (i.e. the vertices of the spaces of operations of $\P$). Let $\ovl{c} := (c_1,\cdots,c_{{m}};c)$ and $\ovl{d} := (d_1,\cdots,d_{{n}};d)$ be two $C$-sequences, and let $\mu \in \P(\ovl{c})$ and $\nu\in\P(\ovl{d})$ be two operations. The data of a morphism $\mu \lrar \nu$ in $\Tw(\P)$ consist of
	
	\smallskip
	
	$\bullet$ a map $f:\langle n \rangle \lrar \langle m \rangle$ in $\Fin_*$,
	
		\smallskip
	
	$\bullet$ a tuple of operations  $\alpha = (\alpha_0, \alpha_1,\cdots,\alpha_m) \in \Map^{f}_{\textbf{Ib}^{\P}}(\ovl{c} , \ovl{d})$, and
	
		\smallskip
	
	$\bullet$ an edge $p : \nu \lrar \alpha^{*}(\mu)$ in $\P(\ovl{d})$ where $\alpha^{*} : \P(\ovl{c}) \lrar \P(\ovl{d})$ is the map corresponding to $\alpha$ via the simplicial functor structure of $\P^{\Inf} : \textbf{Ib}^{\P} \lrar \Set_\Delta$.
\end{rem}

The mapping spaces of $\Tw(\P)$ are given as follows. As above, we regard $\mu \in \P(\ovl{c})$ and $\nu\in\P(\ovl{d})$ as objects of $\Tw(\P)$.

\begin{prop}(\cite{Hoang})\label{p:mappingtwp} {Suppose that $\P$ is fibrant}. There is a canonical homotopy equivalence 
	$$ \{\nu\} \times^{\h}_{\P(\ovl{d})}\Map_{\textbf{Ib}^{\P}}(\ovl{c} , \ovl{d}) \lrarsimeq \Map_{\Tw(\P)}(\mu , \nu) $$
	in which the map $\Map_{\textbf{Ib}^{\P}}(\ovl{c} , \ovl{d}) \lrar \P(\ovl{d})$ is given by the composed map
	$$ \Map_{\textbf{Ib}^{\P}}(\ovl{c} , \ovl{d}) \x{\P^{\Inf}}{\lrar} \Map_{\Set_\Delta}(\mathcal{P}(\ovl{c}) , \mathcal{P}(\ovl{d})) \x{ev_\mu}{\lrar} \mathcal{P}(\ovl{d}) $$
	with $ev_\mu$ being the evaluation at $\mu$.
\end{prop}

The two most basic examples are: $\Tw(\Ass) \simeq \Delta$ and $\Tw(\Com) \simeq \Fin_*^{\op}$. 

\begin{notns} For each $1\leq n < \infty$, we denote by $\E_n$ the simplicial version of the \textbf{little $n$-discs operad}, and denote by $\E_\infty$ a fibrant and  $\Sigma$-cofibrant model for the commutative operad $\Com$.
\end{notns}

As in {\cite{Hoang}}, we have an $\infty$-categorical filtration of the nerve of $\Fin_*^{\op}$:
\begin{equation}\label{eq:nervefin}
	\sN(\Delta) \simeq \Tw(\E_1) \lrar \Tw(\E_2) \lrar \cdots \lrar \Tw(\E_\infty) \simeq \sN(\Fin_*^{\op}).
\end{equation}

{Next, we write $\Sp$ for} the \textbf{$\infty$-category of spectra}, defined as the stabilization of the $\infty$-category of pointed simplicial sets. One of the main interests in the construction $\Tw(-)$ is as follows. 

\begin{thm}(\cite{Hoang})\label{t:maintan} {Suppose that $\P$ is fibrant}. There is an equivalence of $\infty$-categories
	\begin{equation}\label{eq:maintain}
		\T_{\P^{\Inf}}\IbMod(\P)_\infty \simeq \Fun(\Tw(\P) , \Sp).
	\end{equation}
	Consequently, {when $\P$ is in addition $\Sigma$-cofibrant,} there is a chain of equivalences of $\infty$-categories:
	\begin{gather*}
		\T_\P\Op(\Set_\Delta)_\infty \simeq \T_\P\Op_C(\Set_\Delta)_\infty \simeq \T_{\P^{\sB}}\BMod(\P)_\infty \\
		\simeq \T_{\P^{\Inf}}\IbMod(\P)_\infty \simeq \Fun(\Tw(\P) , \Sp).
	\end{gather*}
\end{thm}

\begin{rem}\label{r:acompu} Let $M\in \T_{\P^{\Inf}}\IbMod(\P)$ be a prespectrum whose datum consists of for each $k\geq0$ a sequence of maps $\P^{\Inf} \lrar M_{k,k} \lrar \P^{\Inf}$. We will write $\widetilde{M}:\Tw(\P) \lrar \Sp$ for the functor corresponding to $M$. Unwinding the definitions, $\widetilde{M}$ is given at each operation $\mu\in \P(\ovl{c})$ by taking $\widetilde{M}(\mu)\in\Sp$ to be a prespectrum with 
	$$ \widetilde{M}(\mu)_{k,k} \simeq \{\mu\} \times^{\h}_{\P(\ovl{c})} M_{k,k}(\ovl{c}).$$
	Moreover, for a morphism $\mu \lrar \nu$ in $\Tw(\P)$ (with $\nu\in\P(\ovl{d})$) that lies above a morphism $\ovl{c} \lrar \ovl{d}$ in $\Map_{\textbf{Ib}^{\mathcal{P}}}(\ovl{c} , \ovl{d})$, we get a natural map $\widetilde{M}(\mu)_{k,k} \lrar \widetilde{M}(\nu)_{k,k}$ induced by the commutative square
	$$ \xymatrix{
		M_{k,k}(\ovl{c}) \ar[r]\ar[d] & \P(\ovl{c}) \ar[d] \\
		M_{k,k}(\ovl{d}) \ar[r] & \P(\ovl{d}) \\
	}$$
	where the vertical maps are due to the infinitesimal $\P$-bimodule structures on $M_{k,k}$ and $\P^{\Inf}$. 
\end{rem}

\begin{example}  When $\P=\E_1$ (resp. $\E_\infty$), various tangent categories at this operad are equivalent to $\Fun(\sN(\Delta) , \Sp)$ (resp. $\Fun(\sN(\Fin_*^{\op}) , \Sp)$). In particular, the identification 
	$$ \T_{\E_\infty}\Op(\Set_\Delta)_\infty \simeq \Fun(\sN(\Fin_*^{\op}) , \Sp)$$ 
	proves that the \textbf{stabilization of $\infty$-operads} is equivalent to $\Fun(\sN(\Fin_*^{\op}) , \Sp)$, because $\E_\infty$ is a terminal object in the $\infty$-category $\Op(\Set_\Delta)_\infty$.
\end{example} 

\subsection{Cohomology theories for simplicial operads and algebras over them}\label{s:cotan}

By Theorem \ref{t:maintan}, the Hochschild complex and various cotangent complexes of $\P$ can be represented {as  functors} $\Tw(\P) \lrar \Sp$, provided that $\P$ is fibrant and $\Sigma$-cofibrant. As in {\cite[$\S$6.3]{Hoang}}, we let $\F_\P : \Tw(\P) \lrar \Sp$ denote the image of $\rL_{\P^{\sB}} \in \T_{\P^{\sB}}\BMod(\P)_\infty$ through the identification
$$ \T_{\P^{\sB}}\BMod(\P)_\infty  \simeq \Fun(\Tw(\P) , \Sp) .$$
{Note that} $\F_\P$ coincides with the image of $\P\circ\sS_C^{\bullet}\in\T_{\P^{\Inf}}\IbMod(\P)$ through the equivalence \eqref{eq:maintain}. Alternatively, $\F_\P$ agrees with the image of $\rL_{\P}[1] \in \T_\P\Op(\Set_\Delta)_\infty$ under the identification $$\T_\P\Op(\Set_\Delta)_\infty \simeq \Fun(\Tw(\P) , \Sp)$$
{(cf. $\S$\ref{s:rv})}. In order to give an explicit description of $\F_\P$, we will {need} the following.

\begin{prop}\label{p:destw} Let $f:\P \lrar \Q$ be a map between fibrant and $\Sigma$-cofibrant simplicial operads. Then the derived functor of the right adjoint $$f^* : \T_\Q\Op(\Set_\Delta) \lrar \T_\P\Op(\Set_\Delta)$$ sends $\rL_{\Q}$ to $\rL_{\P}$.
	\begin{proof} It will suffice to show that the derived functor of the right adjoint $$f^* : \T_{\Q^{\Inf}}\IbMod(\Q) \lrar \T_{\P^{\Inf}}\IbMod(\P)$$ sends $\Q\circ\sS_D^{\bullet}$ to $\P\circ\sS_C^{\bullet}$, where $D$ refers to the set of colors of $\Q$. Here by construction the above functor is induced by the right adjoint $$ f^* : \IbMod(\Q)_{\Q^{\Inf}//\Q^{\Inf}} \lrar \IbMod(\P)_{\P^{\Inf}//\P^{\Inf}}$$
		that takes each $M \in \IbMod(\Q)_{\Q^{\Inf}//\Q^{\Inf}}$ to $f^*(M) := \P^{\Inf} \times_{(\Q^{\Inf})^*} M^*$, where the pullback is taken in $\IbMod(\P)$ and $M^*$ denotes the infinitesimal $\P$-bimodule with $$M^*(c_1,\cdots,c_n;c) := M(f(c_1),\cdots,f(c_n);f(c)),$$
		(while the case of $(\Q^{\Inf})^*$ is defined similarly). Unwinding the definitions, we just need to verify for each $m\geq0$ the existence of a natural weak equivalence in $\IbMod(\P)$ of the form
		$$ \P\circ\sS_C^{m} \lrarsimeq \P^{\Inf} \times^{\h}_{(\Q^{\Inf})^*} (\Q\circ\sS_D^{m})^*$$
		(see also {\cite[Corollary 2.4.8]{YonatanBundle}}). Moreover, after having fixed a Kan-model for the simplicial $m$-sphere  $\sS^{m}$, it suffices to show that the canonical map $$\P\circ\sS_C^{m} \lrar \P^{\Inf} \times_{(\Q^{\Inf})^*} (\Q\circ\sS_D^{m})^*$$ is an isomorphism. This is verified because we have on each level a Cartesian square of the form
		$$ \xymatrix{
			\P(c_1,\cdots,c_n;c)\times (\sS^{m})^{\times n} \ar[r]\ar[d] & \Q(f(c_1),\cdots,f(c_n);f(c))\times (\sS^{m})^{\times n} \ar[d] \\
			\P(c_1,\cdots,c_n;c) \ar[r] & \Q(f(c_1),\cdots,f(c_n);f(c)). \\
		}$$
	\end{proof}
\end{prop}

\begin{rem} {Note that the assertion above does not necessarily hold for base categories $\cS$ which are not Cartesian.}
\end{rem}

{When $\P$ is $\Sigma$-cofibrant, we can assume without loss of generality that $\P$ comes equipped with a map of operads $\Phi : \P \lrar \E_\infty$. We then obtain a consequence as follows.}   
\begin{cor}\label{co:destw} {Suppose that $\P$ is fibrant and $\Sigma$-cofibrant}. The functor $\F_\P$ is equivalent to the composed functor
	$$ \Tw(\P) \x{\Tw(\Phi)}{\lrar} \Tw(\E_\infty) \x{\F_{\E_\infty}}{\lrar} \Sp.$$ 
	\begin{proof} {Under the identification of Theorem \ref{t:maintan}, the functor $\Phi^* : \T_{\E_\infty}\Op(\Set_\Delta) \lrar \T_\P\Op(\Set_\Delta)$ corresponds to the functor 
			$$ \Tw(\Phi)^* : \Fun(\Tw(\E_\infty) , \Sp) \lrar \Fun(\Tw(\P) , \Sp) $$
			given by the restriction along $\Tw(\Phi) : \Tw(\P) \lrar \Tw(\E_\infty)$. We complete the proof by combining this with Proposition \ref{p:destw}.}
	\end{proof}
\end{cor}

\begin{rem}\label{r:equidef0} As presented in {\cite[$\S$6.3]{Hoang}}, the functor $$ \F_{\E_\infty} : \sN(\Fin_*^{\op}) \lrar \Sp $$
	is given on objects by $\F_{\E_\infty}(\l m \r) =\mathbb{S}^{\oplus m}$ (i.e. the $m$-fold coproduct of the sphere spectrum); and moreover, for each map $f : \l n \r \lrar \l m \r$ in $\Fin_*$, the structure map 
	$$ \F_{\E_\infty}(f) : \mathbb{S}^{\oplus \{1,\cdots,m\}} \lrar \mathbb{S}^{\oplus \{1,\cdots,n\}}  $$ 
	is defined by, for each $i\in \{1,\cdots,m\}$, copying the {$i$-th} summand to the summands of position $j\in f^{-1}(i)$ when this fiber is nonempty, or collapsing that summand to the zero spectrum otherwise.
\end{rem}

\begin{rem}\label{r:equidef} The functor $\F_{\E_\infty}$ can be equivalently defined by letting $$\F_{\E_\infty}(\l m \r) = [\l m \r, \mathbb{S}]_*$$ where  $[- , -]_*$ refers to the powering of pointed simplicial sets over spectra, and $\l m \r$ is considered as a (discrete) pointed simplicial set with base point $0$. {Accordingly}, we may think of $\F_{\E_\infty}$ as the spectral version of \textbf{Pirashvili's functor} $$t : \Fin_*^{\op} \lrar \Mod(\textbf{k}),$$
	with $\textbf{k}$ being some commutative ring, {defined by}  $t(\l m \r) := [\l m \r, \textbf{k}]_*$, i.e. the $\textbf{k}$-module of based maps $\l m \r \lrar \textbf{k}$ (where $\textbf{k}$ has base point $0_\textbf{k}$). This functor $t$ played a very active role in the author's works, e.g. \cite{Pirash} and {\cite[Chapter 13]{Lodaycyc}}. We shall revisit this remark later.
\end{rem} 

Here is the fundamental theorem concerning the cotangent complex of simplicial operads.
\begin{thm}(\cite{Hoang})\label{t:cotan} {Suppose that $\P$ is fibrant and $\Sigma$-cofibrant}. Then under the equivalence $$\T_{\P^{\sB}}\BMod(\P)_\infty \simeq \Fun(\Tw(\P) , \Sp),$$ the cotangent complex $\rL_{\P^{\sB}} \in \T_{\P^{\sB}}\BMod(\P)_\infty$ is identified with {the functor} $\F_\P : \Tw(\P) \lrar \Sp$ given on objects by sending each operation $\mu\in\P$ of arity $m$  to $\F_\P(\mu) = \mathbb{S}^{\oplus m}$. Consequently, under the equivalence $$\T_\P\Op(\Set_\Delta)_\infty \simeq \Fun(\Tw(\P) , \Sp),$$
	the cotangent complex $\rL_{\P} \in \T_\P\Op(\Set_\Delta)_\infty$ is identified with $\F_\P[-1]$ the desuspension of $\F_\P$.
\end{thm}

\begin{rem}\label{r:immcons} An immediate consequence is as follows. After sending coefficients into $\Fun(\Tw(\P) , \Sp)$, the Quillen cohomology of $\P$ with coefficients in a given functor $\F : \Tw(\P) \lrar \Sp$ is computed by
	$$ \HHQ^\star(\P ; \F) \simeq \Map_{\Fun(\Tw(\P) , \Sp)}(\F_\P[-1], \F).$$
	Moreover, the $n$-th Quillen cohomology group of $\P$ with coefficients in $\F$ is given by
	$$ \HHQ^n(\P ; \F)  \cong \pi_0\Map_{\Fun(\Tw(\P) , \Sp)}(\F_\P, \F[n+1]).$$
\end{rem}

The statement below proves that Quillen cohomology of $\E_\infty$ can be built up from {that} of {the operads $\E_n$}. Let $\F : \Tw(\E_\infty) \lrar \Sp$ be given and let $\F_n$ denote the composed functor 
$$ \Tw(\E_n) \x{\varphi_n}{\lrar} \Tw(\E_\infty) \x{\F}{\lrar} \Sp.$$
\begin{prop} For every integer $k$ there is an exact sequence of abelian groups
	\begin{equation}\label{eq:milnor}
		0 \lrar \underset{n}{\lim}^1 \HHQ^{k-1}(\E_n ; \F_n) \lrar  \HHQ^k(\E_\infty ; \F)  \lrar \underset{n}{\lim} \HHQ^{k}(\E_n ; \F_n) \lrar 0
	\end{equation}
	where $\lim^1$ denotes the first right derived functor of the limit-functor on towers of abelian groups.
	\begin{proof} The filtration \eqref{eq:nervefin} gives rise to a sequential limit of $\infty$-categories:
		$$ \Fun(\Tw(\E_\infty),\Sp) \lrar \cdots \lrar \Fun(\Tw(\E_2),\Sp) \lrar \Fun(\Tw(\E_1),\Sp).$$
		Combined with Corollary \ref{co:destw}, {this yields} a sequential limit of $\infty$-groupoids:
		\begin{gather*}
			\left[\F_{\E_\infty}, \F[k+1]\right]_{\Sp^{\Tw(\E_\infty)}} \lrar \cdots \lrar [\F_{\E_2}, \F_2[k+1]]_{\Sp^{\Tw(\E_2)}} \\
			\lrar [\F_{\E_1}, \F_1[k+1]]_{\Sp^{\Tw(\E_1)}}.
		\end{gather*}
		(Here for brevity we write $[-,-]_\C$ to denote the mapping space in an {$\infty$-category} $\C$). We consider the above spaces as a pointed space with the zero map as its base point. Now, the  \textit{first Milnor exact sequence} (see e.g. {\cite[$\S$9.3]{BK}}) determines an exact sequence of abelian groups
		\begin{gather*}
			0 \lrar \underset{n}{\lim}^1 \pi_1\left[\F_{\E_n}, \F_n[k+1]\right]_{\Sp^{\Tw(\E_n)}} \lrar  \pi_0\left[\F_{\E_\infty}, \F[k+1]\right]_{\Sp^{\Tw(\E_\infty)}} \\
			\lrar \underset{n}{\lim} \pi_0\left[\F_{\E_n}, \F_n[k+1]\right]_{\Sp^{\Tw(\E_n)}} \lrar 0.
		\end{gather*}
		Applying the formulae given in Remark \ref{r:immcons}, we obtain the expected exact sequence \eqref{eq:milnor}.
	\end{proof}
\end{prop}

{We now turn our attention to the Hochschild cohomology of $\P$ and its associated algebras.}

\begin{prop}\label{p:hoch} {Suppose that $\P$ is fibrant.} Under the equivalence $$\T_{\P^{\Inf}}\IbMod(\P)_\infty \simeq \Fun(\Tw(\P) , \Sp),$$ the Hochschild complex $\rL_{\P^{\Inf}}\in\T_{\P^{\Inf}}\IbMod(\P)_\infty$ is identified with $\ovl{\mathbb{S}} : \Tw(\P) \lrar \Sp$, i.e., the constant functor with value $\mathbb{S}\in\Sp$.
	\begin{proof} {By Proposition \ref{p:hhop}}, $\rL_{\P^{\Inf}}$ can be modeled by the suspension spectrum $\P\times\mathbb{S}$ with
		$$ (\P\times\mathbb{S})_{m,m}(c_1,\cdots,c_n;c) = \P(c_1,\cdots,c_n;c) \times \sS^{m}.$$
		We write $\mathcal{H}_\P : \Tw(\P) \lrar \Sp$ for the functor corresponding to $\rL_{\P^{\Inf}}$. As in Proposition \ref{p:destw}, we may show that the {(derived)} right adjoint $\T_{\E_\infty^{\Inf}}\IbMod(\E_\infty) \lrar \T_{\P^{\Inf}}\IbMod(\P)$ sends $\rL_{\E_\infty^{\Inf}}$ to $\rL_{\P^{\Inf}}$. Thus, as in Corollary \ref{co:destw}, we obtain that $\mathcal{H}_\P$ is equivalent to the composed functor $$ \Tw(\P) \lrar \Tw(\E_\infty) \x{\mathcal{H}_{\E_\infty}}{\lrar} \Sp.$$ Therefore, the proof will be completed after showing that $\mathcal{H}_{\E_\infty} : \sN(\Fin_*^{\op}) \lrar \Sp$ is equivalent to the constant functor with value $\mathbb{S}$. For any object $\l n \r$ of $\Fin_*^{\op}$, by construction $\mathcal{H}_{\E_\infty}(\l n \r)\in \Sp$ is the prespectrum with $\mathcal{H}_{\E_\infty}(\l n \r)_{m,m}$ being given by 
		$$ \mathcal{H}_{\E_\infty}(\l n \r)_{m,m} \simeq \{*\}\times^{\h}_{\E_\infty(n)}(\E_\infty\times \, \mathbb{S})_{m,m}(n) \simeq (\E_\infty\times \, \mathbb{S})_{m,m}(n) \simeq \sS^{m}$$
		(see Remark \ref{r:acompu}). Thus we obtain that $\mathcal{H}_{\E_\infty}(\l n \r) \simeq \mathbb{S}$. Clearly the functor structure maps of $\mathcal{H}_{\E_\infty}$ are all the identity on $\mathbb{S}$. We just verified that $\mathcal{H}_{\E_\infty}$ is equivalent to $\ovl{\mathbb{S}}$.
	\end{proof} 
\end{prop}

Due to the above proposition, we may formulate Hochschild cohomology of $\P$ as follows. 

\begin{thm}\label{t:hoch} {Suppose that $\P$ is fibrant.}  The Hochschild cohomology of $\P$ with coefficients in a given functor $\F : \Tw(\P) \lrar \Sp$ is computed by
	$$ \HHH^\star(\P ; \F) \simeq \Map_{\Fun(\Tw(\P) , \Sp)}(\ovl{\mathbb{S}}, \F).$$
	In particular, the $n$-th Hochschild cohomology group of $\P$ with coefficients in $\F$ is computed by
	$$ \HHH^n(\P ; \F) \cong  \pi_0\Map_{\Fun(\Tw(\P) , \Sp)}(\ovl{\mathbb{S}}, \F[n]) \cong  \pi_{-n} \, \lim\F. $$
\end{thm}

Now, let $A$ be a $\P$-algebra classified by a map $\ell_A : \P \lrar \End_A$ of operads. {Recall that there is an equivalence $\T_A\Alg_\P(\cS) \simeq \T_{A^{\me}}\Mod_\P^A$,  provided that $\P$ is $\Sigma$-cofibrant and $A \in \Alg_\P(\Set_\Delta)$ is cofibrant (see Theorem \ref{t:tanalg}). In light of this, the two cohomologies of $A$ will take coefficients in the same category $\T_{A^{\me}}\Mod_\P^A$.}  

\smallskip

Let $N\in\T_{A^{\me}}\Mod_\P^A$ be a levelwise fibrant $\Omega$-spectrum whose datum consists of for each $k\geq0$ a sequence of maps $A^{\me} \lrar N_{k,k} \lrar A^{\me}$ in $\Mod_\P^A$ such that the structure map $N_{k,k} \lrar A^{\me}$ is a fibration, and in addition, {the following squares are all homotopy Cartesian:}
$$ \xymatrix{
	N_{k,k} \ar[r]\ar[d] & A^{\me} \ar[d] \\
	A^{\me} \ar[r] & N_{k+1,k+1}. \\
}$$
\begin{cons}\label{cons:endan} As discussed in Remark \ref{r:endcons}, the right Quillen functor 
	$$ \End^{\st}_{A,-} : \T_{A^{\me}} \Mod_\P^A \lrar \mathcal{T}_{\P^{\Inf}}\IbMod(\P)$$
	sends $N$ to $\End^{\st}_{A,N}\in \mathcal{T}_{\P^{\Inf}}\IbMod(\P)$ which is a levelwise fibrant $\Omega$-spectrum such that $(\End^{\st}_{A,N})_{k,k}$ is described via a Cartesian square of the form
	$$ \xymatrix{
		(\End^{\st}_{A,N})_{k,k} \ar[r]\ar[d] & \End_{A,N_{k,k}} \ar[d] \\
		\P \ar[r]^{\ell_A} & \End_A. \\
	}$$
	We will make use of Remark \ref{r:acompu} to describe the corresponding functor 
	$$\widetilde{\End}^{\st}_{A,N} : \Tw(\P) \lrar \Sp.$$
	Let $\ovl{c}:=(c_1,\cdots,c_n;c)$ be a $C$-sequence. $\widetilde{\End}^{\st}_{A,N}$ is given on objects by sending each operation $\mu\in\P(\ovl{c})$ to an $\Omega$-spectrum $\widetilde{\End}^{\st}_{A,N}(\mu)\in\Sp$ with 
	$$\widetilde{\End}^{\st}_{A,N}(\mu)_{k,k} = \{\mu\}\times_{\P(\ovl{c})}(\End^{\st}_{A,N})_{k,k}(\ovl{c}) \cong \{\ell_A(\mu)\}\times_{\End_A(\ovl{c})}\End_{A,N_{k,k}}(\ovl{c}).$$
	More explicitly, $\widetilde{\End}^{\st}_{A,N}(\mu)_{k,k}$ agrees with the fiber over $\ell_A(\mu)$ of the simplicial map
	$$ \Map_{\Set_\Delta}(A(c_1)\times\cdots\times A(c_n) , N_{k,k}(c)) \lrar \Map_{\Set_\Delta}(A(c_1)\times\cdots\times A(c_n) , A(c))$$
	that is induced by the structure map $N_{k,k}(c) \lrar  A(c)$ of $N$.
\end{cons}

By having a description of $\widetilde{\End}^{\st}_{A,N}$ as above, and by combining Theorems \ref{t:cotan}, \ref{t:hoch} with the comparison theorems of $\S$\ref{s:op2}, the two cohomologies of $A$ are formulated as follows.

\begin{thm}\label{t:endan} {Suppose that $\P\in\Op_C(\Set_\Delta)$ is $\Sigma$-cofibrant, that $A \in \Alg_\P(\Set_\Delta)$ is cofibrant, and that a levelwise fibrant $\Omega$-spectrum $N\in\T_{A^{\me}}\Mod_\P^A$ is given.}
	
	\smallskip
	
	(i) The Quillen cohomology of $A$ with coefficients in $N$ is computed by
	\begin{gather*}
		 \HHQ^\star(A ; N) \simeq \HHQ^\star(\P^{\sB} ; \End^{\st}_{A,N}) \simeq \Omega\HHQ^\star(\P ; \End^{\st}_{A,N}) \\
		 \simeq \Map_{\Fun(\Tw(\P) , \Sp)}(\F_\P, \widetilde{\End}^{\st}_{A,N}).
	\end{gather*}

	(ii) The Hochschild cohomology of $A$ with coefficients in $N$ is computed by
	$$ \HHH^\star(A ; N) \simeq \HHH^\star(\P ; \End^{\st}_{A,N}) \simeq \Map_{\Fun(\Tw(\P) , \Sp)}(\ovl{\mathbb{S}}, \widetilde{\End}^{\st}_{A,N}).$$
	In particular, the $n$-th Hochschild cohomology group of $A$ with coefficients in $N$ is given by
	$$   \HHH^n(A ; N) \cong \pi_{-n} \, \lim\widetilde{\End}^{\st}_{A,N}.$$
\end{thm}

\begin{example} {Suppose we are given a functor $T : \Fin_*^{\op} \lrar \Mod({\textbf{k}})$ (also referred to as a \textbf{right $\Gamma$-module}), where $\textbf{k}$ is a commutative ring. The $k$-th \textbf{stable cohomotopy group} of $T$ is 
		$$ \pi^{k}T := \Ext_{\Fin_*^{\op}}^{k}(t,T) \cong \pi_0\Map^{\h}_{\Fun(\Fin_*^{\op},\C(\textbf{k}))}(t,T[k])$$
		where $t : \Fin_*^{\op} \lrar \Mod({\textbf{k}})$ refers to Pirashvili's functor (see Remark \ref{r:equidef}). Recall from \cite{Hoang} that there is a natural isomorphism
		$$  \pi^{k}T \cong \HHQ^{k-1} (\E_\infty ; \widetilde{T})$$
		between the $k$-th stable cohomotopy group of $T$ and the $(k-1)$-th Quillen cohomology group of $\E_\infty$ with coefficients in the induced functor $\widetilde{T} : \sN(\Fin_*^{\op}) \lrar \Sp$. This demonstrates that Robinson's \textbf{gamma cohomology} (\cite{Robinson1}) can be naturally encompassed by Quillen cohomology of $\E_\infty$, namely when the coefficient-functors come from right $\Gamma$-modules. Furthermore, combined with the main result of \cite{Robinson}, it shows that the obstruction to an $\E_\infty$-structure is governed by Quillen cohomology of the operad $\E_\infty$ itself.} 
\end{example}

More examples concerning the obtained results will be given in the remainder of the paper.

\subsection{Quillen principle for $\E_n$-algebras}\label{s:quiprin}

{We begin by introducing some terminology and notations.}

\begin{define}\label{d:adeop} {A simplicial operad is \textbf{adequate} if it is fibrant and its spaces of nullary and unary (i.e., $1$-ary) operations are weakly contractible.}
\end{define}

{Typical examples of adequate simplicial operads include the operads $\E_n$ for $0\leq n \leq \infty$.}

\begin{notns}\label{no:xT} {For any object $x$ in an $\infty$-category $\C$ and any spectrum object $T\in\Sp$, we denote by} $x_![T] : \C \lrar \Sp$ (resp. $x_*[T] : \C \lrar \Sp$) the left (resp. right) Kan extension of the embedding $\{T\} \lrar \Sp$ along the inclusion $\{x\} \lrar \C$. 
\end{notns}

\begin{rem}\label{r:xT} Unwinding the definitions, the functor $x_![T] : \C \lrar \Sp$ sends an object $y\in\C$ to $x_![T](y) \simeq \Map_\C(x,y)\otimes T$ where $-\otimes-$ refers to the copowering of simplicial sets over spectra. On other hand, the functor $x_*[T] : \C \lrar \Sp$ sends $y$ to $x_*[T](y) \simeq [\Map_\C(y,x), T]$ where $[-,-]$ denotes the powering of simplicial sets over spectra. In particular, when $T = \mathbb{S}$ then we have $$x_![\mathbb{S}](y) \simeq \Map_\C(x,y)\otimes\mathbb{S} \simeq \Sigma^\infty_+\Map_\C(x,y),$$ and while $$x_*[\mathbb{S}](y) \simeq [\Map_\C(y,x), \mathbb{S}] \simeq (\Sigma^\infty_+\Map_\C(y,x))^\vee,$$
	i.e. the \textit{dual spectrum} of the suspension spectrum $\Sigma^\infty_+\Map_\C(y,x)$.
\end{rem}

\begin{rem}\label{r:xxT} Suppose given a functor $\F : \C \lrar \Sp$. By adjunction, a morphism $x_![T] \lrar \F$ (resp. $\F \lrar x_*[T]$) in $\Fun(\C,\Sp)$ can be represented by a map $T \lrar \F(x)$  (resp. $\F(x) \lrar T$) of spectra. 
\end{rem}

{In this subsection, we will let $\P\in\Op_*(\Set_\Delta)$ be a fixed single-colored simplicial operad.} As usual, we write $\id \in \P(1)$ for the identity operation, and let $\mu_0 \in \P(0)$ be a nullary operation selected arbitrarily. We will regard the two as objects of $\Tw(\P)$.

\begin{rem}\label{r:mu0} When $\P$ is \textit{adequate} in the sense of Definition \ref{d:adeop}, $\mu_0$ represents a terminal object in $\Tw(\P)$ (cf. {\cite[$\S$6.2]{Hoang}}).
\end{rem}

\begin{lem}\label{l:twid} {Suppose that $\P\in\Op_*(\Set_\Delta)$ is adequate}. For any operation $\mu\in\P(m)$ considered as an object of $\Tw(\P)$, there is a weak equivalence of spaces
	$$ \Map_{\Tw(\P)}(\mu,\id) \simeq \P(2) \sqcup \underline{m}$$
	in which $\underline{m} = \{1,\cdots,m\}$ represents a discrete space of cardinality $m$. 
	\begin{proof} Firstly, since $\P$ has a single color, we may identify objects of $\textbf{Ib}^{\mathcal{P}}$ with those of $\Fin_*$. We have in fact a chain of weak equivalences 
		$$ \Map_{\Tw(\P)}(\mu,\id) \simeq \{\id\} \times^{\h}_{\P(1)}\Map_{\textbf{Ib}^{\P}}(\l m \r , \l 1 \r) \simeq \Map_{\textbf{Ib}^{\P}}(\l m \r , \l 1 \r) \simeq \P(2) \sqcup \underline{m}$$ 
		in which the first weak equivalence is due to Proposition \ref{p:mappingtwp}, while the second follows from the assumption that $\P(1)  \simeq *$; and moreover, the last weak equivalence is interpreted as follows. We have by construction that
		$$ \Map_{\textbf{Ib}^{\P}}(\l m \r , \l 1 \r) = \underset{f}{\bigsqcup} \, \P(k_0) \times \P(k_1) \times \cdots \times  \P(k_m)$$
		where the coproduct ranges over those maps $f\in\Hom_{\Fin_*}(\l 1 \r,\l m \r)$ and each $k_i\in\NN$ refers to the cardinality of the fiber $f^{-1}(i)$. Observe now that the summand that corresponds to the constant map $\l 1 \r \lrar \l m \r$ with value $0$ is given by $\P(2)\times \P(0)^{\times m}$, whilst for every $j\in\{1,\cdots,m\}$ the summand corresponding to the map $$[\xi^j : \l 1 \r \lrar \l m \r, 1 \mapsto j]$$ is simply given by $\P(1)\times\P(1) \times \P(0)^{\times (m-1)}$. We complete the proof using the assumption that $\P(0) \simeq \P(1)  \simeq *$.
	\end{proof}
\end{lem}

{The cotangent complex of $\P$ can be described through its space of binary operations, as follows.}

\begin{prop}\label{p:fundp} {Suppose that $\P$ is adequate}. Then there is a fiber sequence in $\Fun(\Tw(\P),\Sp)$ of the form
	\begin{equation}\label{eq:fundp}
		\F_\P \x{\iota}{\lrar} \id_*[\mathbb{S}] \x{\rho}{\lrar} \ovl{(\Sigma^\infty_+\P(2))^\vee} 
	\end{equation}
	where $\ovl{(\Sigma^\infty_+\P(2))^\vee}$ refers to the constant functor $\Tw(\P) \lrar \Sp$ with value $(\Sigma^\infty_+\P(2))^\vee$ (see Remark \ref{r:xT} for notation).
	\begin{proof} Concretely, the map $\F_\P \x{\iota}{\lrar} \id_*[\mathbb{S}]$ is classified by the identity $\mathbb{S} = \F_\P(\id) \lrar \mathbb{S}$, and the map $\rho$ is given as follows. Observe first that a map $\id_*[\mathbb{S}] \lrar \ovl{(\Sigma^\infty_+\P(2))^\vee}$ is completely determined by a map in $\Sp$ of the form
		$$ \id_*[\mathbb{S}](\mu_0) \lrar (\Sigma^\infty_+\P(2))^\vee,$$
		due to the fact that $\mu_0$ is a terminal object in $\Tw(\P)$ (see Remark \ref{r:mu0}). Combining Remark \ref{r:xT} with Lemma \ref{l:twid}, we obtain for each operation $\mu\in\P(m)$ a chain of equivalences in $\Sp$: 
		$$ \id_*[\mathbb{S}](\mu) \simeq [\Map_{\Tw(\P)}(\mu,\id), \mathbb{S}] \simeq [\P(2) \sqcup \underline{m}, \mathbb{S}] \simeq (\Sigma^\infty_+\P(2))^\vee\oplus \mathbb{S}^{\oplus m}.$$
		In particular, we obtain that $\id_*[\mathbb{S}](\mu_0) \simeq (\Sigma^\infty_+\P(2))^\vee$. Accordingly, we would take $\rho$ to be the map classified by the identity on $ (\Sigma^\infty_+\P(2))^\vee$. Now the sequence \eqref{eq:fundp} is given at each $\mu\in\P(m)$ by the sequence
		$$ \xymatrix{
			\F_\P(\mu) \ar[r]\ar[d]_{=} & \id_*[\mathbb{S}](\mu) \ar[r]\ar[d]^{\simeq} & \ovl{(\Sigma^\infty_+\P(2))^\vee} \, (\mu) \ar[d]^{=} \\
			\mathbb{S}^{\oplus m} \ar[r] & (\Sigma^\infty_+\P(2))^\vee \oplus \mathbb{S}^{\oplus m} \ar[r] & (\Sigma^\infty_+\P(2))^\vee , \\
		} $$
		which is clearly a fiber sequence in $\Sp$. The proof is therefore completed.
	\end{proof}
\end{prop}

\begin{example} {Recall that $\Tw(\E_\infty) \simeq \sN(\Fin_*^{\op})$, and hence, in this situation the fiber sequence \eqref{eq:fundp} takes the form}
	\begin{equation}\label{eq:fundcom}
		\F_{\E_\infty} \x{\iota}{\lrar} \l 1 \r_*[\mathbb{S}] \lrar  \ovl{\mathbb{S}}.
	\end{equation}
	{More explicitly, the map $\iota$ is given at each object $\l m \r$ by the embedding $$\F_{\E_\infty}(\l m \r) \simeq \mathbb{S}^{\oplus m}  \lrar \mathbb{S}^{\oplus (m+1)} \simeq \l 1 \r_*[\mathbb{S}](\l m \r)$$
		whose image leaves out the summand $\mathbb{S}$ corresponding to the constant map $\l 1 \r \rar \l m \r$ with value $0$. This is the spectral version of an assertion obtained by Pirashvili {\cite[$\S$1.4]{Pirash}} (see also Remark \ref{r:equidef}).}  
\end{example}

\begin{rem}\label{r:fundp} {From the fiber sequence} \eqref{eq:fundp}, we obtain for each functor $\F : \Tw(\P) \lrar \Sp$ a fiber sequence of mapping spaces:
	\begin{gather*}
		\Map_{\Sp} \left( (\Sigma^\infty_+\P(2))^\vee,\lim\F  \right) \lrar \Map_{\Fun(\Tw(\P),\Sp)}(\id_*[\mathbb{S}],\F) \\ \lrar \Map_{\Fun(\Tw(\P),\Sp)}(\F_\P,\F).
	\end{gather*}
	{By Theorem \ref{t:hoch} and Remark \ref{r:immcons}, one may view this sequence as establishing a relation between Quillen and Hochschild cohomologies of operads in general.} 
\end{rem}

{Our main interest is the connection between the two cohomologies of $\E_n$-spaces. First, let us recall the steps to establish the fiber sequence \eqref{eq:Qprin1} from \cite{Hoang}, providing additional clarification.} 

\smallskip

As above, we write $\id\in\E_n(1)$ for the identity operation, and write $\mu_0\in\E_n(0)$ for the unique nullary operation of $\E_n$. Let us consider the functor $$(\mu_0)_![\mathbb{S}] : \Tw(\E_n) \lrar \Sp.$$
We let $\mu \in \E_n(k)$ be a fixed $k$-ary operation regarded as an object of $\Tw(\E_n)$. Due to Proposition \ref{p:mappingtwp}, we obtain that 
$$ \Map_{\Tw(\E_n)}(\mu_0,\mu) \simeq \{\mu\} \times^{\h}_{\E_n(k)}\E_n(k+1)$$
where the map $\E_n(k+1) \lrar \E_n(k)$ is given by $\E_n(k+1)\ni\nu \mapsto \nu\circ_1 \mu_0$ (i.e. given by deleting the first $n$-disc in $\nu$ when considering $\nu$ as a rectilinear embedding of $k+1$ disjoint $n$-discs in another $n$-disc). Thus we obtain further a weak equivalence $$\Map_{\Tw(\E_n)}(\mu_0,\mu) \simeq \underset{k}{\bigvee} \sS^{n-1}$$
where the latter space is the $k$-fold wedge sum of the $(n-1$)-sphere. Therefore, as in {\cite[$\S$6.3]{Hoang}} we have that
$$ (\mu_0)_![\mathbb{S}](\mu) \, \simeq \, \Map_{\Tw(\E_n)}(\mu_0,\mu) \otimes \mathbb{S} \, \simeq \, (\underset{k}{\bigvee} \sS^{n-1}) \otimes \mathbb{S} \, \simeq \, \mathbb{S} \oplus \left(\mathbb{S}[n-1]\right)^{\oplus k}.$$
If we consider the functor $(\mu_0)_!\left[\mathbb{S}[-n+1]\right] : \Tw(\E_n) \lrar \Sp$ instead, then we get that
$$ (\mu_0)_!\left[\mathbb{S}[-n+1]\right](\mu) \, \simeq \, \Omega^{n-1} \left[ (\mu_0)_![\mathbb{S}](\mu) \right] \, \simeq \, \mathbb{S}[-n+1] \oplus \mathbb{S}^{\oplus k}.$$

\smallskip

We also regard the functor $\id_*[\mathbb{S}] : \Tw(\E_n) \lrar \Sp$. According to Lemma \ref{l:twid}, and along with the fact that $\E_n(2)\simeq\sS^{n-1}$, we obtain that
\begin{gather*}
	\id_*[\mathbb{S}](\mu) \, \simeq \, [\Map_{\Tw(\E_n)}(\mu,\id), \mathbb{S}] \, \simeq \, [\sS^{n-1} \sqcup \, \underline{k} \, , \mathbb{S}] \, \\
	 \simeq \,  [\sS^{n-1}, \mathbb{S}] \oplus [\underline{k} \, , \mathbb{S}]  \, \simeq \, \left(\mathbb{S}\oplus\mathbb{S}[-n+1]\right) \oplus \mathbb{S}^{\oplus k}.
\end{gather*}

Let us now consider the canonical morphism in $\Fun(\Tw(\E_n),\Sp)$: $$\varepsilon : (\mu_0)_!\left[\mathbb{S}[-n+1]\right] \lrar \id_*[\mathbb{S}]$$ that is classified by the inclusion $\ovl{\varepsilon} : \mathbb{S}[-n+1] \lrar \mathbb{S}\oplus\mathbb{S}[-n+1] \simeq \id_*[\mathbb{S}](\mu_0)$. More explicitly, the morphism  $\varepsilon$ is given at the operation $\mu$ by an embedding of the form
$$ (\mu_0)_!\left[\mathbb{S}[-n+1]\right](\mu) \, \simeq \, \mathbb{S}[-n+1] \oplus \mathbb{S}^{\oplus k} \lrar \mathbb{S}\oplus\mathbb{S}[-n+1] \oplus \mathbb{S}^{\oplus k} \, \simeq \, \id_*[\mathbb{S}](\mu)  $$
whose image leaves out the distinguished summand $``\mathbb{S}$'' that is contained in $$[\sS^{n-1} , \mathbb{S}] \simeq \mathbb{S}\oplus\mathbb{S}[-n+1].$$

{From the above analysis, we deduce the following.}
\begin{prop}\label{p:fundEn} There is a fiber sequence in $\Fun(\Tw(\E_n),\Sp)$ of the form
	\begin{equation}\label{eq:fundEn}
		(\mu_0)_!\left[\mathbb{S}[-n+1]\right] \x{\varepsilon}{\lrar} \id_*[\mathbb{S}]  \x{\pi}{\lrar} \ovl{\mathbb{S}}
	\end{equation}
	where as usual $\ovl{\mathbb{S}}$ denotes the constant functor $\Tw(\E_n) \lrar \Sp$ with value $\mathbb{S}$, and the morphism $\pi$ is canonically classified by the projection $$\id_*[\mathbb{S}](\mu_0) \simeq \mathbb{S}\oplus\mathbb{S}[-n+1] \lrar \mathbb{S}.$$
\end{prop}

\begin{example} {We recall that $\Tw(\E_1) \simeq \sN(\Delta)$. Therefore, in this case we get a fiber sequence in $\Fun(\sN(\Delta),\Sp)$ of the form}
	$$ [0]_![\mathbb{S}] \x{\varepsilon}{\lrar} [1]_*[\mathbb{S}]  \x{\pi}{\lrar} \ovl{\mathbb{S}}.$$
	Let us give an explicit description of the morphism $\varepsilon$ above, which is an interesting one. By definition, $\varepsilon$ is classified by the inclusion $\ovl{\varepsilon} :   \mathbb{S} \lrar \mathbb{S}\oplus\mathbb{S} = [1]_*[\mathbb{S}]([0])$. We will let $\ovl{\varepsilon}(\mathbb{S})$ be the summand that corresponds to the map $\beta^1 : [0] \lrar [1] , 0 \mapsto 1$. Now the morphism $\varepsilon$ is given at each object $[k]$ of $\sN(\Delta)$ by an embedding of the form 
	$$ \varepsilon([k]) : \Hom_\Delta([0],[k]) \otimes \mathbb{S} \, = \, \mathbb{S}^{\oplus (k+1)} \lrar  \mathbb{S}^{\oplus (k+2)} \, = \, [\Hom_\Delta([k],[1]), \mathbb{S}].$$ 
	For each map $[\alpha^i : [0] \lrar [k], 0 \mapsto i]$, we write $(\alpha^i, \mathbb{S})$ for $\mathbb{S}$ itself regarded as a summand of the left hand side indexed by $\alpha^i$. Moreover, for a map $f : [k] \lrar [1]$, we write $\mathbb{S}^{f}$ to denote $\mathbb{S}$ in its role as a summand of the right hand side that corresponds to $f$. Accordingly, $\varepsilon([k])$ is determined by copying each $(\alpha^i, \mathbb{S})$ to those summands $\mathbb{S}^{f}$ with $f$ being such that $f\circ\alpha^i = \beta^1$ (i.e. such that $f(i)=1$). In particular, the image of $\varepsilon([k])$ leaves out only the summand $\mathbb{S}^{\const_0}$ in which  $\const_0$ signifies the constant map $[k] \lrar [1]$ with value $0$.
\end{example}

Applying the fiber sequence \eqref{eq:fundp} {for} $\P = \E_n$, we obtain a fiber sequence in $\Fun(\Tw(\E_n), \Sp)$ of the form
$$ \F_{\E_n} \lrar \id_*[\mathbb{S}] \lrar \ovl{(\Sigma^\infty_+\E_n(2))^\vee} \simeq \ovl{\mathbb{S}}\oplus\ovl{\mathbb{S}}[-n+1].$$
Combining this with the fiber sequence \eqref{eq:fundEn}, we obtain a corollary as follows.
\begin{cor}\label{co:fundEn} There is a diagram of Cartesian squares in $\Fun(\Tw(\E_n),\Sp)$ of the form
	$$ \xymatrix{
		\F_{\E_n} \ar[r]\ar[d] & (\mu_0)_!\left[\mathbb{S}[-n+1]\right] \ar[r]\ar[d] & \id_*[\mathbb{S}] \ar[d] \\
		0 \ar[r] & \ovl{\mathbb{S}}[-n+1] \ar[r] & \ovl{\mathbb{S}}\oplus\ovl{\mathbb{S}}[-n+1] \\
	} $$
	where the last vertical map is induced by the identity $$\id_*[\mathbb{S}](\mu_0) \simeq \mathbb{S}\oplus\mathbb{S}[-n+1] \lrar \mathbb{S}\oplus\mathbb{S}[-n+1].$$
\end{cor}

In particular, the left square of the above diagram is also Cartesian. It allows us to recover a result presented in {\cite[$\S$6.3]{Hoang}} as follows.

\begin{thm}(\cite{Hoang}) There is  a fiber sequence in $\Fun(\Tw(\E_n),\Sp)$ of the form
	\begin{equation}\label{eq:Qprin1}
		(\mu_0)_![\mathbb{S}] \lrar \ovl{\mathbb{S}} \lrar \F_{\E_n}[n].
	\end{equation}
\end{thm}

As we have seen, the functor $\F_{\E_n}$ represents the cotangent complex of $\E_n$, and while $\ovl{\mathbb{S}}$ represents the Hochschild complex of $\E_n$ (see Proposition \ref{p:hoch}). In light of this, we will refer to the fiber sequence  \eqref{eq:Qprin1} as the \textbf{Quillen principle for the operad} $\E_n$.

\smallskip

{We now explain} how the Quillen principle for $\E_n$-algebras can be formulated, {relying on the following lemma. For a single-colored simplicial operad $\P$,} we write $\eta_\P : \Free_\P^{\Inf}(\mathcal{E}_*) \lrar \P^{\Inf}$ for the map in $\IbMod(\P)$ classified by {a given nullary operation} $\mu_0 \in \P(0)$ (see $\S$\ref{s:not} for notations). 

\begin{lem}\label{l:preQprin} {Suppose that $\P$ is fibrant}. Under the equivalence of $\infty$-categories
	$$  \Fun(\Tw(\P),\Sp)  \simeq \T_{\P^{\Inf}}\IbMod(\P)_\infty,$$
	the map $(\mu_0)_![\mathbb{S}]\lrar \ovl{\mathbb{S}}$ is identified with the canonical map $\Sigma^\infty_+(\eta_\P) \lrar \rL_{\P^{\Inf}}$ in which as usual $\Sigma^\infty_+(\eta_\P)$ denotes the image of $\eta_\P$ through the functor $$\Sigma^\infty_+ : \IbMod(\P)_{/\P^{\Inf}} \lrar \T_{\P^{\Inf}}\IbMod(\P).$$
	\begin{proof} We let $\textbf{S} := (\Set_\Delta)_\infty$ denote the $\infty$-category of spaces. We also write $(\mu_0)_! : \textbf{S} \lrar \Fun(\Tw(\P),\textbf{S})$ for the induction along the inclusion $\{\mu_0\} \lrar \Tw(\P)$. We have a commutative diagram of the form
		$$ \xymatrix{
			\textbf{S} \ar[r]^{(\mu_0)_! \;\;\;\;\;\;\;\;\;\;\;\;\;\;}\ar[d]_{\simeq} & \Fun(\Tw(\P),\textbf{S}) \ar[r]^{\Sigma^\infty_+\;\,}\ar[d]^{\simeq} & \Fun(\Tw(\P),\Sp) \ar[d]^{\simeq} \\
			\Sigma_*(\textbf{S})_{/\mathcal{E}_*} \ar[r]_{\Free_\P^{\Inf} \;\;\;\;\;\;\;\;} & {\IbMod(\P)_\infty}_{/\P^{\Inf}} \ar[r]_{\Sigma^\infty_+\;\;} & \T_{\P^{\Inf}}\IbMod(\P)_\infty \\
		} $$
		in which the left vertical functor is the obvious equivalence sending each space $X$ to the $\Sigma_*$-object given by $X$ itself concentrated in level $0$ and equipped with an evident map to $\mathcal{E}_*$, and while the middle vertical functor is the equivalence induced by the straightening-unstraightening constructions (cf. {\cite[$\S$6.3]{Hoang}}). Let us start with the terminal object $\Delta^0 \in \textbf{S}$. Clearly the image of $\Delta^0$ through the composed horizontal functor is identified with $(\mu_0)_![\mathbb{S}] \in \Fun(\Tw(\P),\Sp)$. On the other hand, observe that the left vertical functor sends $\Delta^0$ to nothing but $\Id_{\mathcal{E}_*} \in \Sigma_*(\textbf{S})_{/\mathcal{E}_*}$, whose image through $\Free_\P^{\Inf}$ is exactly the map $\eta_\P$. Now, the commutativity of the diagram proves that $(\mu_0)_![\mathbb{S}]$ is identified with $\Sigma^\infty_+(\eta_\P)$ under the equivalence $\Fun(\Tw(\P),\Sp)  \simeq \T_{\P^{\Inf}}\IbMod(\P)_\infty$.
		
		\smallskip
		
		Moreover, by Proposition \ref{p:hoch}, we obtained that $\ovl{\mathbb{S}}$ is identified with the Hochschild complex $\rL_{\P^{\Inf}}$. To complete the proof, we argue that the two maps $(\mu_0)_![\mathbb{S}]\lrar \ovl{\mathbb{S}}$ and $\Sigma^\infty_+(\eta_\P) \lrar \rL_{\P^{\Inf}}$ are both induced by the operation $\mu_0$.
	\end{proof}
\end{lem}

Combining the above lemma with the fiber sequence \eqref{eq:Qprin1}, and along with Theorem \ref{t:cotan} we obtain a corollary as follows.

\begin{cor}\label{co:Qprin2} There is a fiber sequence in $\T_{\E_n^{\Inf}}\IbMod(\E_n)$ of the form
	\begin{equation}\label{eq:Qprin2}
		\Sigma^\infty_+(\eta_{\E_n}) \lrar \rL_{\E_n^{\Inf}}  \lrar \rL_{\E_n}[n+1].
	\end{equation}
	Here we use the same notation for the derived image of $\rL_{\E_n}\in\T_{\E_n}\Op(\Set_\Delta)$ under the equivalence $\T_{\E_n}\Op(\Set_\Delta) \simeq \T_{\E_n^{\Inf}}\IbMod(\E_n)$.
\end{cor}

\begin{rem} The fiber sequence \eqref{eq:Qprin2} is in fact the operadic source of \eqref{eq:Qprin1}, and in particular, can be viewed as the standard version of the Quillen principle for $\E_n$. One might think of another proof of the existence of \eqref{eq:Qprin2} that is more direct; nevertheless, it may be much more complicated.
\end{rem}

We are now in position to formulate the \textbf{Quillen principle for $\E_n$-algebras}. Let $A$ be an  $\E_n$-algebra. We let $\Free_A(*)$ denote the free $A$-module over $\E_n$ generated by a singleton, and let $\eta_A : \Free_A(*) \lrar A^{\me}$ be the map in $\Mod_{\E_n}^A$ classified by the unit of $A$. As usual, we have a canonical map $\Sigma^\infty_+(\eta_A) \lrar \rL_{A^{\me}}$ where $\Sigma^\infty_+(\eta_A)$ refers to the image of $\eta_A$ through the functor $$\Sigma^\infty_+ : (\Mod_{\E_n}^A)_{/A^{\me}} \lrar \T_{A^{\me}}\Mod_{\E_n}^A.$$
\begin{thm}\label{t:Qprin3} {Suppose that $A\in\Alg_{\E_n}(\Set_\Delta)$ is cofibrant.} There is a fiber sequence in $\T_{A^{\me}}\Mod_{\E_n}^A$ of the form
	\begin{equation}\label{eq:Qprin3}
		\Sigma^\infty_+(\eta_A) \lrar \rL_{A^{\me}} \lrar \rL_A[n]
	\end{equation}
	where we use the same notation for the derived image of $\rL_A\in\T_{A}\Alg_{\E_n}(\Set_\Delta)$ under the equivalence $\T_{A}\Alg_{\E_n}(\Set_\Delta) \simeq \T_{A^{\me}}\Mod_{\E_n}^A$. Equivalently, $\rL_A[n]$ is weakly equivalent to $\rL_{A^{\me}/\Free_A(*)}$ the relative cotangent complex of $\eta_A$.
	\begin{proof} First, note that the fiber sequence \eqref{eq:Qprin2} can be equivalently rewritten as 
		$$ \Sigma^\infty_+(\eta_{\E_n}) \lrar \rL_{\E_n^{\Inf}}  \lrar \rL_{\E_n^{\sB}}[n] $$ 
		where we use the same notation for the derived image of $\rL_{\E_n^{\sB}}\in\T_{\E_n^{\sB}}\BMod(\E_n)$ under the equivalence $\T_{\E_n^{\sB}}\BMod(\E_n) \simeq \T_{\E_n^{\Inf}}\IbMod(\E_n)$. Let us now apply the left Quillen functor $$(-)\circ^{\st}_{\E_n} A : \mathcal{T}_{\E_n^{\Inf}}\IbMod(\E_n)  \lrar \T_{A^{\me}}\Mod_{\E_n}^A$$ (see Notation \ref{no:circend}) to that fiber sequence. Then we may complete the proof by using Corollary \ref{r:freeA}, which proves that $ \Sigma^\infty_+(\eta_{\E_n})\circ^{\st}_{\E_n} A \simeq \Sigma^\infty_+(\eta_A)$; and along with Proposition \ref{p:bialg}, from which we obtain the following weak equivalences: $ \rL_{\E_n^{\Inf}}\circ^{\st}_{\E_n} A  \simeq \rL_{A^{\me}}$    and     $\rL_{\E_n^{\sB}}[n] \circ^{\st}_{\E_n} A \simeq \rL_A[n]$.
	\end{proof}
\end{thm}

\begin{cor}\label{co:Qprin3} {Suppose that $A\in\Alg_{\E_n}(\Set_\Delta)$ is cofibrant and} suppose given a fibrant object $\N \in \T_{A^{\me}}\Mod_{\E_n}^A$. Then there is a fiber sequence relating the two cohomologies of $A$ of the form
	$$ \Omega^n\HHQ^\star(A;\N) \lrar \HHH^\star(A;\N) \lrar |\N|$$
	in which $|\N| := \{1_A\} \times^{\h}_A \N_{0,0}$, i.e. the homotopy fiber over $1_A$ of the structure map $\N_{0,0} \lrar A^{\me}$.
	\begin{proof} We have in fact a commutative diagram of the form
		\begin{equation*}
			\xymatrix{
				\Map^{\h}(\rL_A[n],\N) \ar[r]\ar[d]_{\simeq} & \Map^{\h}(\rL_{A^{\me}},\N) \ar[r]\ar[d]^{\simeq} & \Map^{\h}(\Sigma^\infty_+(\eta_A),\N) \ar[d]^{\simeq} \\
				\Omega^n\HHQ^\star(A;\N) \ar[r] & \HHH^\star(A;\N) \ar[r] & |\N| \\
			} 
		\end{equation*}
		where the derived mapping spaces are taken in $\T_{A^{\me}}\Mod_{\E_n}^A$, and the first two weak equivalences follow immediately from the definitions. To verify the third weak equivalence, we argue as follows. First, by adjunction we have a weak equivalence
		$$ \Map^{\h}_{\T_{A^{\me}}\Mod_{\E_n}^A}(\Sigma^\infty_+(\eta_A),\N) \simeq  \Map^{\h}_{(\Mod_{\E_n}^A)_{/A^{\me}}}(\eta_A, \Omega^\infty_+\N).$$
		Next, the space on the right is weakly equivalent to the homotopy fiber over $\eta_A$ of the map
		$$ \Map^{\h}_{\Mod_{\E_n}^A}(\Free_A(*),\N_{0,0}) \lrar \Map^{\h}_{\Mod_{\E_n}^A}(\Free_A(*),A^{\me}).$$
		Moreover, clearly this homotopy fiber is weakly equivalent to $|\N|$.
	\end{proof}
\end{cor}

{Next} we provide a version of the Quillen principle for $\E_n$-algebras in a stable category. Suppose we are given a base category $\cS$ that is sufficient in the sense of {\cite[$\S$3]{Hoang}}, and equipped with a \textbf{weak monoidal Quillen adjunction} $\adjunction*{}{\Set_\Delta}{\cS}{}$ (cf. {\cite[$\S$3.2]{Schwede}}). We will use the same notation for the $\cS$-enriched version of the operad $\E_n$. 

\begin{thm}\label{t:Qprin56} Let $A \in \Alg_{\E_n}(\cS)$ be a cofibrant $\E_n$-algebra in $\cS$. There is a fiber sequence in $\T_{A^{\me}}\Mod_{\E_n}^A(\cS)$ of the form
	\begin{equation}\label{eq:Qprin5}
		\Sigma^\infty_+(\eta_A) \lrar \rL_{A^{\me}} \lrar \rL_A[n]
	\end{equation}
	(i.e. the same as in Theorem \ref{t:Qprin3}). When $\cS$ is in addition strictly pointed stable, then there is a fiber sequence in $\Mod_{\E_n}^A(\cS)$ of the form
	\begin{equation}\label{eq:Qprin6}
		\Free_A(*) \x{\eta_A}{\lrar} A^{\me} \lrar \rL_A[n].
	\end{equation}
	\begin{proof} Clearly both the cotangent and Hochschild complexes of $\E_n$ are preserved through the changing-bases functor $\Set_\Delta \lrar \cS$. It implies that the same fiber sequence as in Corollary \ref{co:Qprin2} holds for the $\cS$-enriched version of $\E_n$. Then the same argument as in the proof of Theorem \ref{t:Qprin3} can be applied to obtain the fiber sequence \eqref{eq:Qprin5}. 
		
		\smallskip
		
		Now assume further that $\cS$ is stable. In particular, the composed functor $$\ker \circ \, \Omega^\infty : \T_{A^{\me}}\Mod_{\E_n}^A(\cS) \lrarsimeq \Mod_{\E_n}^A(\cS)$$ is a right Quillen equivalence (see Theorem \ref{t:tanalg}). Then we may obtain the fiber sequence \eqref{eq:Qprin6} by using {\cite[Corollary 2.2.4]{Yonatan}}, which demonstrates that the derived image of the map $\Sigma^\infty_+(\eta_A) \lrar \rL_{A^{\me}}$ through the functor $\ker \circ \, \Omega^\infty$ is nothing but the map $\Free_A(*) \x{\eta_A}{\lrar} A^{\me}$.
	\end{proof}
\end{thm}

\begin{examples}\label{ex:Qprin56} {The most basic example is the case of \textbf{symmetric spectra}}. Another noteworthy example is the case $\cS = \C(\textbf{k})$ the category of dg modules over some commutative ring $\textbf{k}$. Nonetheless, it is a bit subtle that $\C(\textbf{k})$ {does not come equipped with a lax monoidal left Quillen functor from} $\Set_\Delta$. Instead, we have a composed functor 
	$$ \Set_\Delta \x{\textbf{k}[-]}{\lrar} \sMod_{\textbf{k}} \x{\sN}{\lrar} \C_{\geq0}(\textbf{k}) \x{\iota}{\lrar} \C(\textbf{k})$$
	in which $\textbf{k}[-]$ refers to the free $\textbf{k}$-module functor, while $\sN$ denotes the normalization functor, and $\iota$ is the obvious embedding. This composed functor takes $\E_n$ to its \textbf{differential graded (dg) version}. In the above composition, while $\textbf{k}[-]$ and $\iota$ are left adjoints of weak monoidal Quillen adjunctions,  $\sN$ is the only functor that does not satisfy the expected property. Despite this, the latter induces a right Quillen equivalence between operadic infinitesimal bimodules,  as presented in \cite{Hoang1}. Due to this, we may show that an analogue of \eqref{eq:Qprin2} holds for the dg version of $\E_n$. Then we may deduce the existence of \eqref{eq:Qprin6} for \textbf{dg $\E_n$-algebras}, as in Theorem \ref{t:Qprin56}.
\end{examples}

\newpage

\bibliographystyle{amsplain}

\begin{thebibliography}{10}
	
	
		\smallskip
		
	
	\bibitem{Turchin} G. Arone and V. Turchin, {\it On the rational homology of high-dimensional analogues of spaces of long knots}, Geom. Topol. 18(3) (2014) 1261-1322.
	
	\smallskip
	
	\bibitem{White} M. Batanin and D. White, {\it Left Bousfield localization without left properness}, {J. Pure
		Appl. Algebra. 228(6)} (2024), DOI: 10.1016/j.jpaa.2023.107570.
	
	\smallskip
	
	\bibitem{BauWi} H. J. Baues and G. Wirsching, {\it Cohomology of small categories}, J. Pure
	Appl. Algebra. 38 (1985) 187-211.			
	
	\smallskip
	
	\bibitem{BergerMoerdijk} C. Berger and I. Moerdijk, {\it On the derived category of an algebra over an operad},  Georgian Math. J. 16(1) (2009) 13-28.
	
	
	\smallskip
	
	\bibitem{BK} A. K. Bousfield and D. M. Kan, {\it Homotopy limits, completions and localizations}, Lecture Notes in Mathematics. 304, Springer-Verlag, Berlin-New York (1972).
	
	\smallskip
	
	\bibitem{Caviglia} G. Caviglia, {\it A model structure for enriched coloured operads}, {preprint, 	arXiv:1401.6983 (2014)}.
	
	\smallskip
	
	
	\bibitem{Julien} J. Ducoulombier, B. Fresse and V. Turchin, {\it Projective and Reedy model category structures for (infinitesimal) bimodules over and operad}, Appl. Categ. Struct. 30 (2022) 825-920.
	
	\smallskip
	
	
	\bibitem{Francis} J. Francis, {\it The tangent complex and Hochschild cohomology of $\sE_n$-rings}, Compos. Math. 149(3) (2013) 430-480.
	
	\smallskip
	
	\bibitem{Fresse1} B. Fresse, {\it Modules over operads and functors}, Lecture Notes in Mathematics. 1967, Springer-Verlag, Berlin (2009).
	
	\smallskip
	
	
	\bibitem{Guti}  J. J. Gutiérrez and R. M. Vogt, {\it A model structure for coloured operads in symmetric spectra}, Math. Z. 270 (2012) 223-239.
	
	\smallskip
	
	\bibitem{YonatanCotangent} Y. Harpaz, J. Nuiten and M. Prasma, {\it The abstract cotangent complex and Quillen cohomology of enriched categories}, J. Topol. 11(3) (2018) 752-798.	
	
	\smallskip				
	
	
	\bibitem{Yonatan} Y. Harpaz, J. Nuiten and M. Prasma, {\it  Tangent categories of algebras over operads}, Isr. J. Math. 234(2) (2019) 691-742.
	
	\smallskip
	
	
	\bibitem{YonatanBundle} Y. Harpaz, J. Nuiten and M. Prasma, {\it The tangent bundle of a model category}, Theory Appl.
	Categ. 34 (2019) 1039-1072.
	
	
	\smallskip
	
	
	
	
	
	
	
	
	{\bibitem{Hoang}      T. Hoang, {\it Quillen cohomology of enriched operads}, Adv. Math. 465
		(2025),  DOI: 10.1016/j.aim.2025.110151.}
	
	\smallskip
	
	
	
	\bibitem{Hoang1}      T. Hoang, {\it Operadic Dold-Kan correspondence and mapping spaces between enriched operads}, preprint, arXiv:2312.07906 (2023).
	
	\smallskip
	
	\bibitem{Hoch} G. Hochschild, {\it On the cohomology groups of an associative algebra}, Ann. of Math., Second Series, 46 (1) (1945) 58-67.
	
	\smallskip
	
	
	\bibitem{Hovey} M. Hovey, {\it Model categories}, Mathematical Surveys and Monographs. 63, Amer. Math. Soc. (1999).
	
	
	\smallskip
	
	
	
	\bibitem{Lodaycyc} J.-L. Loday, {\it Cyclic homology}, Second Edition, Grundlehren Math.Wiss. 301, Springer, Berlin-Heidelberg (1998).		
	
	\smallskip
	
	
	\bibitem{Loday} J.-L. Loday and B. Vallette, {\it Algebraic Operads}, Grundlehren Math.Wiss. 346, Springer, Heidelberg  (2012).
	
	
	\smallskip
	
	
	
	\bibitem{Lurieha} J. Lurie, {\it Higher algebra}, preprint, available at {author’s homepage \url{http://www.math.ias.edu/~lurie/}} (2011).
	
	\smallskip
	
	\bibitem{Vallette} S. Merkulov and B. Vallette, {\it Deformation theory of representations of prop(erad)s. II}, J. Reine Angew. Math. 636 (2009) 1-52.
	
	\smallskip
	
	
	\bibitem{Pirash}  T. Pirashvili, {\it Hodge decomposition for higher order Hochschild Homology}, Ann. Sci. Ecole Norm. Sup. 33(2) (2000) 151-179.
	
	
	\smallskip
	
	
	
	\bibitem{Rich} T. Pirashvili and B. Richter, {\it Hochschild and cyclic homology via
		functor homology}, K-theory. 25(1) (2002) 39-49.
	
	\smallskip
	
	
	\bibitem{Quillen} D. Quillen, {\it Homotopical algebra}, Lecture Notes in Mathematics. 43, Springer-Verlag (1967).
	
	\smallskip
	
	
	
	\bibitem{Rezk} C. Rezk, {\it Spaces of algebra structures and cohomology of operads}, PhD Thesis, Massachusetts Institute of Technology (1996).
	
	\smallskip
	
	\bibitem{Robinson1} A. Robinson, {\it Gamma homology, Lie representations and $\E_\infty$ multiplications}, Invent. Math. 152(2) (2003) 331-348.
	
	\smallskip
	
	
	{\bibitem{Robinson} A. Robinson, {\it $E_\infty$-obstruction theory}, Homol. Homotopy Appl. 20(1) (2018) 155-184.}
	
	\smallskip
	
	
	\bibitem{Schwede} S. Schwede and B. Shipley, {\it Equivalences of monoidal model categories}, Algebr. Geom. Topol. 3(1) (2003) 287-334.
	
	\smallskip
	
	\bibitem{Spitzweck} M. Spitzweck, {\it Operads, algebras and modules in general model categories}, preprint, arXiv:math/0101102 (2001).
	
	\smallskip
	
	
	\bibitem{Toen} B. Toën, {\it The homotopy theory of dg-categories and derived Morita theory}, Invent. Math. 167(3) (2007)  615-667.
	
	
		




\end{thebibliography}

\end{document}